%% file: main.tex
\documentclass[11pt]{article}
\usepackage{amsthm,amssymb,amsmath}
\usepackage{graphicx}
\usepackage[english]{babel}
\usepackage{enumerate}
\usepackage{url}
\usepackage{color}
\usepackage{subfigure}
\usepackage{authblk}
\usepackage{hyperref}
\usepackage{cleveref}
\usepackage{tikz}

\DeclareMathOperator    \aff                {aff}

\newcommand{\bd}                {\partial}
\DeclareMathOperator    \cl                  {cl}
\DeclareMathOperator    \conv            {conv}

\DeclareMathOperator    \intr               {int}

\DeclareMathOperator    \proj              {proj}

\DeclareMathOperator    \relint             {rel\,int}

\DeclareMathOperator    \vol                 {vol}

\newcommand{\R}{\mathbb R}
\newcommand{\Q}{\mathbb Q}
\newcommand{\Z}{\mathbb Z}

\renewcommand{\epsilon}{\varepsilon}

\newtheorem{theorem}{Theorem}
\newtheorem{corollary}[theorem]{Corollary}
\newtheorem{lemma}[theorem]{Lemma}
\newtheorem{proposition}[theorem]{Proposition}
\newtheorem{remark}[theorem]{Remark}

\newtheorem{conjecture}[theorem]{Conjecture}
\newtheorem{example}[theorem]{Example}

\newtheorem{assumption}[theorem]{Assumption}

\newtheorem{definition}[theorem]{Definition}

\newtheoremstyle{oraclestyle}%
  {3pt}%
  {3pt}%
  {}%
  {}%
  {\bfseries}%
  {:}%
  {.5em}%
  {\thmname{#1}\thmnumber{ #2}\thmnote{ (#3)}}%

\theoremstyle{oraclestyle}
\newtheorem{oracle}{Oracle}
\crefname{oracle}{oracle}{oracles}
\Crefname{oracle}{Oracle}{Oracles}

\def\ve#1{\mathchoice{\mbox{\boldmath$\displaystyle\bf#1$}}
{\mbox{\boldmath$\textstyle\bf#1$}}
{\mbox{\boldmath$\scriptstyle\bf#1$}}
{\mbox{\boldmath$\scriptscriptstyle\bf#1$}}}

\renewcommand{\a}{{\ve a}}
\renewcommand{\b}{{\ve b}}
\renewcommand{\c}{{\ve c}}
\newcommand{\w}{{\ve w}}
\newcommand{\x}{{\ve x}}
\newcommand{\y}{{\ve y}}
\newcommand{\z}{{\ve z}}

\usepackage{pict2e,picture}
\newsavebox\CBox
\newlength\CLength
\def\Circled#1{\sbox\CBox{#1}%
  \ifdim\wd\CBox>\ht\CBox \CLength=\wd\CBox\else\CLength=\ht\CBox\fi
    \makebox[1.2\CLength]{\makebox(0,1.2\CLength){\put(0,0){\circle{1.2\CLength}}}%
    \makebox(0,1.2\CLength){\put(-.5\wd\CBox,0){#1}}}}

\usepackage{xspace}
\newcommand{\HC}{Boundary Hyperplane Cover\xspace}
\newcommand{\CCPsub}{convex/concave polyhedral subdivision}
\newcommand{\CCPsubs}{convex/concave polyhedral subdivisions }

\newcommand{\SEP}[2]{\texttt{SEP}_{#1,\,#2}}
\newcommand{\SEPint}[2]{\texttt{SEP}^{\,\circ}_{#1,\,#2}}
\newcommand{\FEAS}[2]{\texttt{FEAS}_{#1,\,#2}}
\newcommand{\WEAK}[2]{\texttt{WeakSep}_{#1,\,#2}}

\newcommand{\PC}{\mathcal P\setminus \mathcal C}

\newcommand{\IE}{\texttt{{F}}}
\usepackage{algorithm}
\usepackage{algorithmic}

\title{Complexity of Integer Programming in Reverse Convex Sets via \HC%
\thanks{R. Hildebrand was partially funded by  AFOSR grant FA9550-21-1-0107. Any opinions, findings, and conclusions or recommendations expressed in this material are those of the authors and do not necessarily reflect the views of the Air Force Office of Scientific Research. A. G{\"o}{\ss} thanks the Deutsche Forschungsgemeinschaft for their support within project A05 in the \emph{Sonderforschungsbereich/Transregio 154 Mathematical Modelling, Simulation and Optimization using the Example of Gas Networks}.}}
\author[1]{Robert Hildebrand} %
\author[2]{Adrian G{\"o}{\ss}}

\affil[1]{Grado Department of Industrial and Systems Engineering, Virginia Tech}
\affil[2]{Analytics \& Optimization Lab, University of Technology Nuremberg (UTN)}

\usepackage[backend=biber, sorting=none, maxnames=999]{biblatex}

\AtEveryBibitem{
    \ifboolexpr{
        not test {\iffieldundef{doi}}
    }
    {\clearfield{url} \clearfield{eprint}}
    {}
}

\usepackage{mathrsfs}
\renewcommand{\H}{H}%

\begin{document}

\maketitle

\begin{abstract}
We study the complexity of identifying the integer feasibility of reverse convex sets. We present various settings where the complexity can be either NP-Hard or efficiently solvable when the dimension is fixed.  Of particular interest is the case of bounded reverse convex constraints with a polyhedral domain.  
We introduce a structure, \emph{\HC}, that permits this problem to be solved in polynomial time in the number of removed sets, hyperplanes, and encoding size provided the dimension is fixed.
\end{abstract}

\newpage 
\section{Introduction}

Our work is motivated by analyzing special structures to help solve Integer Nonlinear Programs.  We consider convex sets $K_1, \dots, K_l$ and $C_1, \dots, C_m$, subsets of $\R^{n}$.  We will study sets of the form
\begin{equation}
\label{eq:main}
    \Big((K_1\cup \dots \cup K_l) \setminus (C_1 \cup \dots \cup C_m)\Big) \cap \Z^n.
\end{equation}
We will show how non-convexities in this formulation can be isolated in various cases.  We will begin with simpler examples, and then work towards a more robust characterization of sets where this operation is possible.

We are primarily interested in the case of \emph{reverse convex programming}.  In this setting, we consider non-convex sets of the form 
\begin{equation*}
    \{ \x \in P : f_i(\x) \geq 0, \text{ for } i =1, \dots, m\},
\end{equation*}
where $P = \{\x \in \R^n : A\x \leq \b\}$ is a rational polyhedron and the sets $\{\x \in \R^n : f_i(\x) \geq 0\}$ are \emph{reverse convex sets}, i.e., the complement sets $C_i = \{ \x \in \R^n : f_i(\x) \leq 0\}$ are convex.

We propose a structure called \emph{\HC} that allows us to decompose the nonlinearities in the solution space to a collection of polyhedra interacting solely with a single non-convex element.  Specifically, a \HC for a collection of convex sets $C_1, \dots, C_m$ is a set of hyperplanes whose union contains all pairwise intersections of the boundaries of the convex sets.  See Definition~\ref{defn:bhc} for a detailed formal definition.    This property arises naturally in a variety of contexts, but may require a special structure in order to occur.  For example, the intersection of any two balls (radii may vary) can be separated conveniently.  However, this property does not follow for any two general ellipsoids.

\setcounter{footnote}{-1}

\begin{figure}[h]
    \centering
    \subfigure[]{
\includegraphics[scale=0.15, trim=300pt 100pt 500pt 220pt, clip]{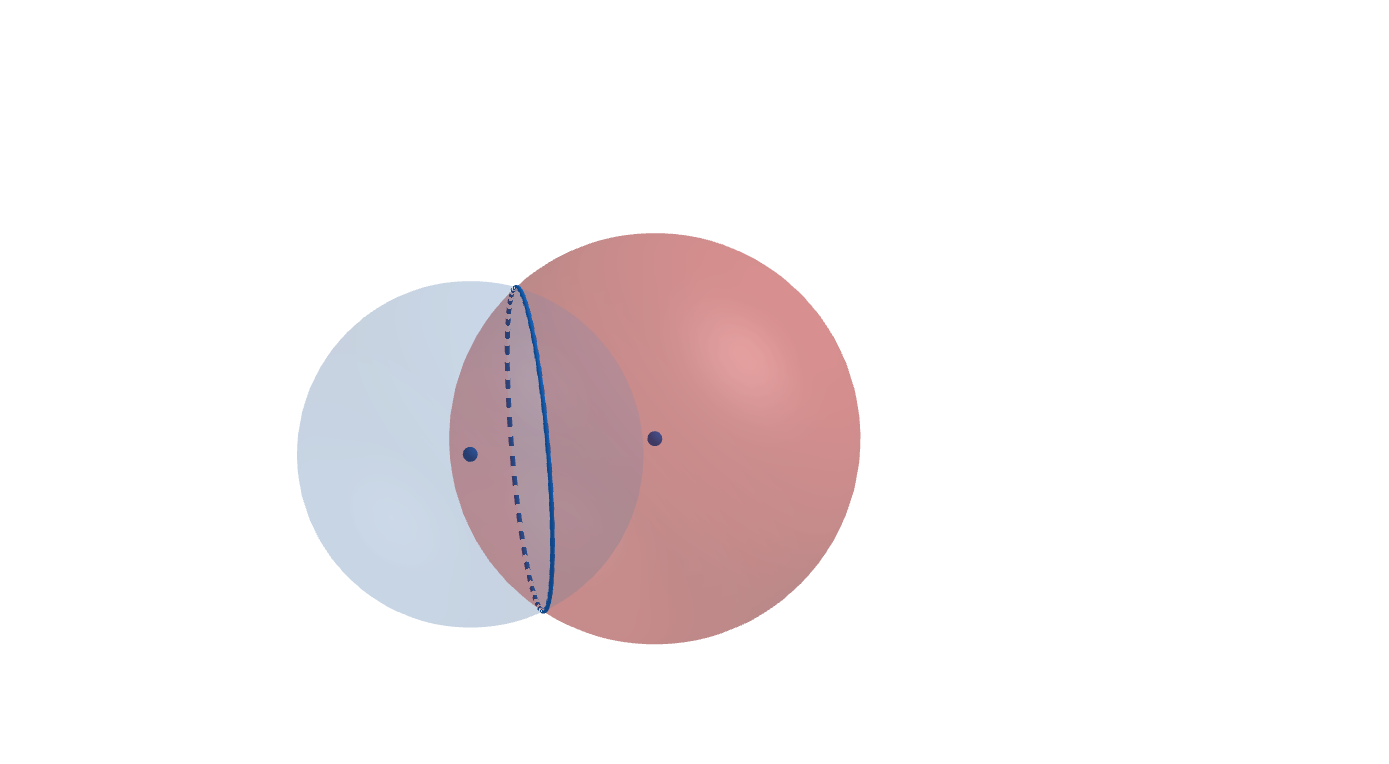}

    }
    \subfigure[]{
        \includegraphics[scale=0.09, trim=500pt 300pt 1300pt 700pt, clip]{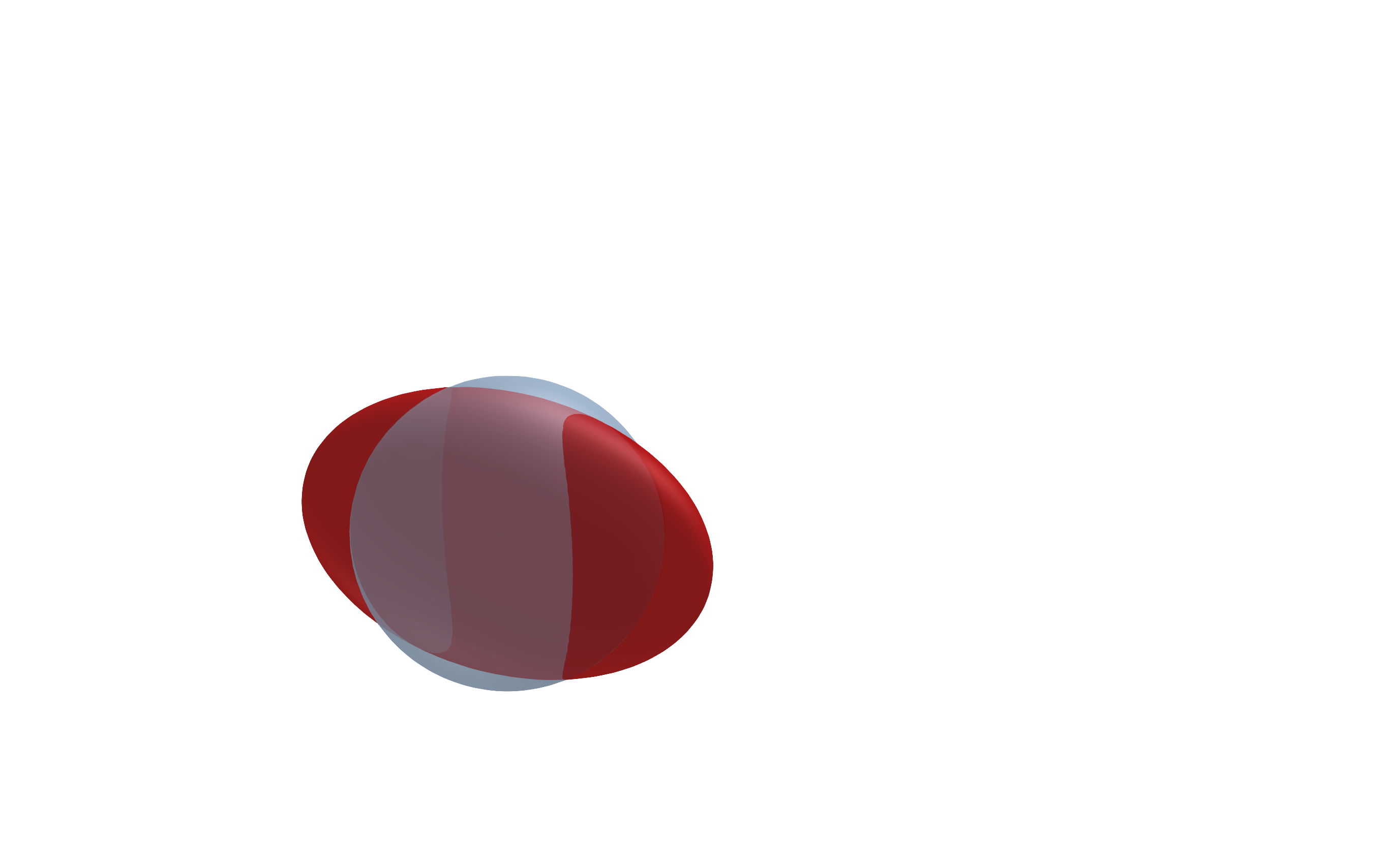}
    }
    \subfigure[]{
        \includegraphics[scale=0.09]{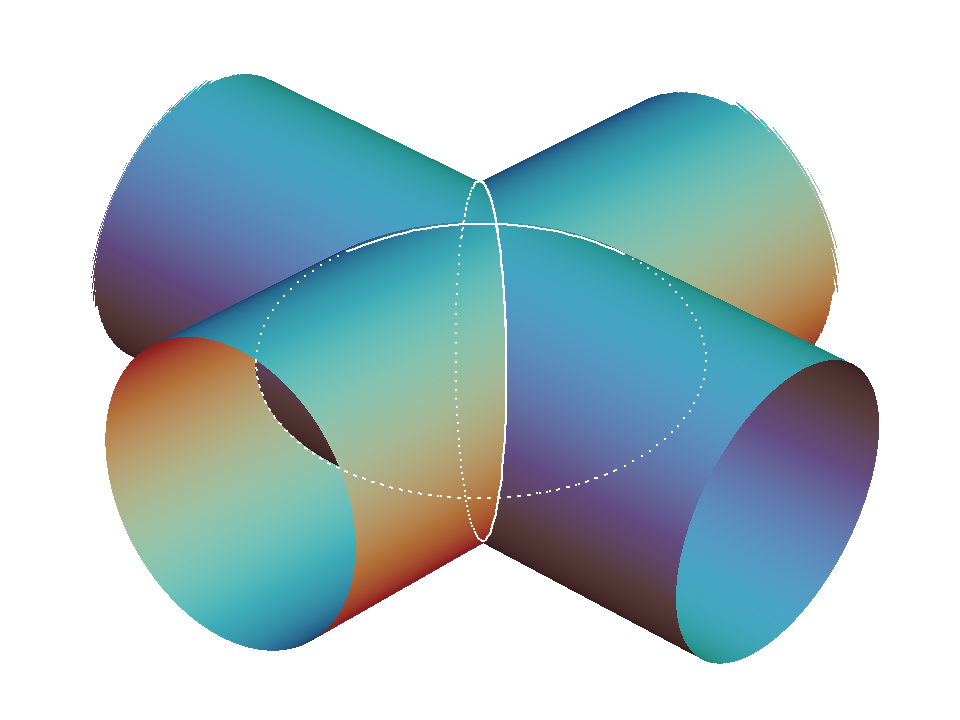}
    }
    \subfigure[]{
        \includegraphics[scale=0.09]{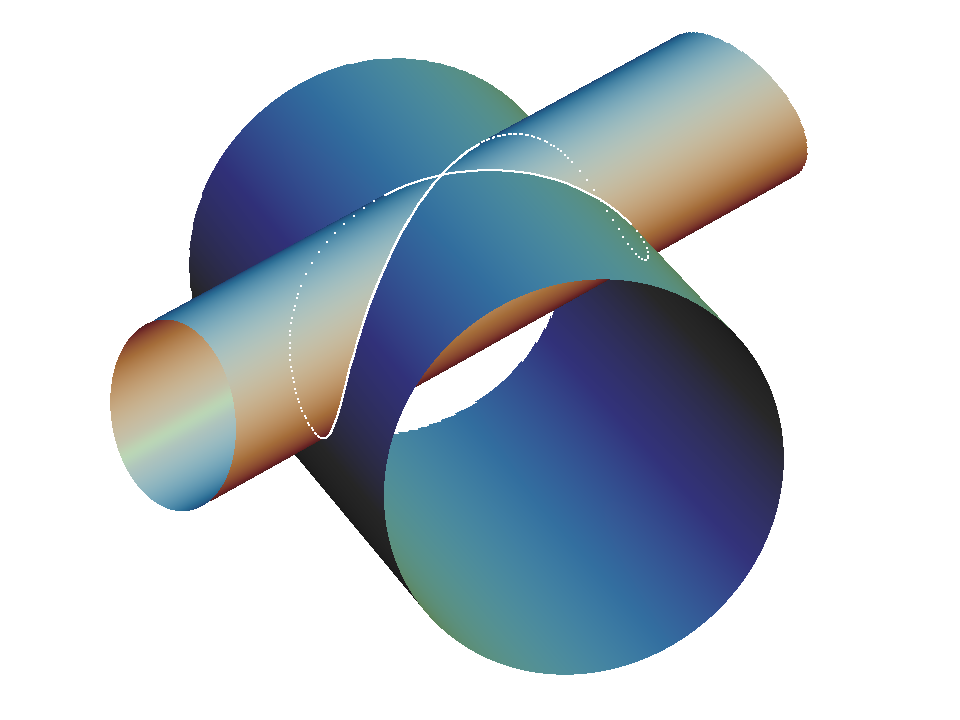}
    }
    \caption[(a) The intersection of two balls with different radii. (b) The intersection of two ellipsoids where the \HC property does not exist. (c) Cylinders that can be separated using two hyperplanes. (d) A case where this doesn't happen.]{(a) The intersection of two balls with different radii. (b) The intersection of two ellipsoids where the \HC property does not exist. (c) Cylinders that can be separated using two hyperplanes. (d)  A case where this property fails.\footnotemark}
    \label{fig:subfigures}
\end{figure}

Cylinders are another nice (related) example; in this case, two hyperplanes are needed in the separation.  \Cref{fig:subfigures} shows examples where this does and doesn't hold.

\footnotetext{Images of ellipsoids were created using \emph{Geogebra}~\cite{geogebra_ellipsoids}. 
The cylinders were modified from~\cite{virtualmathmuseum_cylinders}.}

\HC allows us to apply operations to these subproblems that can only be applied in special contexts.  This property then serves as a feature to define tractable integer nonlinear programs. We apply this idea primarily to the complexity of integer feasibility in reverse convex sets.

We give a brief literature review of these areas, followed by an outline of the paper.

\subsection{Literature: Complexity of integer feasibility}
We distinguish between convex and non-convex sets here.  We will focus on integer feasibility, although much of this literature extends to mixed-integer settings.  See~\cite{oertel, dadush2012integer, Basu2022}.

\paragraph{Convex Integer Feasibility}
Lenstra's seminal work~\cite{lenstra_integer_1983} showed that mixed-integer linear programming can be solved in time $f(n) \mathrm{poly}( \phi)$ where $n$ is the number of integer variables, $f$ is a function depending only on $n$, and $\phi$ is the binary encoding size of the problem (objective and constraints).  There have been a number of works improving the complexity and extending this work to other settings.  In particular, Lenstra's novel approach extends quite readily to the convex integer programming setting~\cite{Grtschel1993}.  There have been generalizations to polynomial and semialgebraic settings~
\cite{Heinz2005, HildebrandK13, khackiyan_integer_2000} and significant work by Dadush on the general oracle settings with optimizing $f(n)$~\cite{dadush2012integer}.  Recent work uses Dadush's framework to achieve an even greater speedup~\cite{rothvoss-subspace}. Other recent works: The work in~\cite{alberto-lenstra} details complexity when there are quadratic constraints, while~\cite{Basu2022} studies the information complexity of the problem.  Near-optimal information complexities are obtained using a center point technique initially developed in~\cite{oertel}.

Alternatively, there is a highly distinct direction based on Graver Bases; see~\cite{graver} and references therein.  This line of work achieves a parameterized complexity based on parameters such as treewidth - a measure of sparsity related to the problem based on variable interactions in the constraints. These tools apply in the linear and also the separable convex case (but with linear constraints), and typically in the pure integer setting.

\paragraph{Non-convex integer feasibility}
There is relatively little work in this area. One tool is the ability to compute the integer or mixed-integer hull of a polyhedron~\cite{cook-hartmann-kannan-mcdiarmid-1992, Hildebrand2015}.  For problems of minimizing a concave function over the integer or mixed-integer points in a polyhedron, one can find an optimal solution at a vertex of the mixed-integer hull.  Thus, those techniques can be applied.  This approach was utilized in  
\cite{delpia_integer_2013} and \cite{cubic-integer-plane} to develop the first-ever algorithms for minimizing quadratic and cubic functions over integer points in polyhedra in the plane (only 2 variables).  
In general dimensions, there is work to minimize non-convex quadratics over integer points in polyhedra provided that there are certain bounds on the size of numbers used in the description of the problem~\cite{zemmer-thesis, Lokshtanov15}.  See also~\cite{Eiben2019} for some related parameterized complexity of the problem.

\subsection{Literature: Reverse Convex Sets}
There are a number of applications that have reverse convex sets (see, e.g., ~\cite{Jacobsen2009}). One such recent work studies the obnoxious facility location problem~\cite{Kalczynski2021}.  
Given a set of communities as vectors in the two-dimensional plane, it aims for placing a certain number of facilities in the same plane.
Those facilities are noisy or have pollution, both of which are transferred via air. 
Therefore, one seeks to establish a minimal distance in the sense of the (squared) Euclidean norm between each facility and community.
In addition, each pair of facilities are required to have a minimal distance to each other.

Both described sets of constraints are reverse convex ones, as they lower bound the squared Euclidean norm of some linear term in the variables by a constant.
If, in addition, the facilities are required to be placed on certain sites, which can be represented as discrete points, an additional integral component is included.
This results in a mixed-integer reverse convex problem.

On some higher level, there exist other areas of mathematical programming which directly or indirectly handle reverse convex constraints. 
For instance, in a difference of convex functions (DC) program, the occurring constraints are indeed differences of convex functions, as the name indicates. 
Each such constraint can be rewritten in terms of a convex and a reverse convex constraint by introduction of one additional variable, leading to a reverse convex program.
We refer the interested reader to~\cite{tuy1986dc} for a more detailed explanation in an approach for DC programming.
The recent work of~\cite{leudtke-reverse-convex} also summarizes the mentioned procedure as a motivation but explicitly studies cutting planes in the context of reverse convex sets.

With regards to the complexity of mixed-integer reverse convex programming, there is a line of work identifying critical constraints and also solving a related partitioning problem~\cite{OBUCHOWSKA201258,OBUCHOWSKA2014129}.

In this paper, we stick to the pure integer setting.  Many of the techniques will adapt to the mixed-integer setting, but, as noted in~\cite{OERTEL2014424}, details about mixed-integer problems require some discussion of approximate oracles to be precise.  We aim to handle such settings in future work.

\subsection{Contributions}
We provide the first-ever complexity results for reverse convex integer programming in dimensions $n \geq 3$ when there are $m\geq 2$ reverse convex constraints.  We establish that, with some additional structure,  we can obtain polynomial-time complexities when the dimension is fixed.  We also study the appearance of this structure and the flexibility of these settings.  

Our approach develops a polynomial-time method to decompose the feasibility problem into simpler ones, which requires the separation of integer points in different convex sets.  As a by-product of developing the tools for this, we present statements about the separation and optimization over the integer points on the interiors of convex sets. 
To the best of our knowledge, this has not been stated before.

\subsection{Outline}
In~\Cref{sec:preliminaries} we lay the groundwork for the paper by providing notation, assumptions and oracles associated with respective justifications, as well as useful statements to keep the main results clear.
In particular, oracles for separation and feasibility are presented and distinguished by their field of operation, i.e., the integers or the reals.

\Cref{sec:hard-simple-cases} provides hardness results on examples when considering problems even in dimension two that suffer from unspecific bounds. 
The section also features the simple example of removing polytopes, which provides a first glimpse into the main idea of this article.
The idea is detailed in~\Cref{sec:main-continuous} and proposes a tractable decomposition of the sets to exclude based on a \emph{\HC}. 
This decomposition relies on (exact) weak separation.
As this operation can cause numerical troubles in practice, a variant based on integer separation is given consecutively in~\Cref{sec:main-integer}, using a similar yet numerically valid approach.
Specifically, we discuss separation of integer points over the interior of convex sets in order to provide the decomposition.
This also involves new feasibility and optimization statements over interiors, which have never been stated before.
All results obtained allow to draw a statement about the computational complexity of structured reverse convex integer programming (\Cref{sec:complexity-rev-convex}).

Lastly, we also provide a short analysis of sets that show the structure of a \HC in~\Cref{sec:characterizations}, before concluding with future work in this regard.
In the appendix, we provide alternative or complementary proofs of several statements.

\section{Preliminaries}
\label{sec:preliminaries}

\subsection{Notation}
For $S \subseteq \R^n$, let $\intr(S)$ denote the interior of $S$, let $\relint(S)$ denote the relative interior of $S$, and let $\bd S$ denote the boundary of $S$.  That is, $\bd S = \cl(S) \setminus \relint(S)$ where $\cl(S)$ is the topological closure of $S$.

For $\epsilon > 0$, let $$
B(S,\epsilon):= \{\x \in \R^n : \|\x - \y\|_2 \leq \epsilon \text{ for some } \y \in S\}.
$$
For $\bar \x \in \R^n$, let $B(\bar \x,\epsilon) := B(\{\bar \x\}, \epsilon)$.  This is exactly the Euclidean ball of radius $
\epsilon$ centered at $\bar \x$.
If we instead consider the $\infty$-norm, we denote it with $B_\infty(\bar \x, \varepsilon)$.

For any $\a \in \R^n \setminus \{{\bf 0}\}$, $b \in \R$, denote the associated hyperplane and halfspaces as
$$
\H_{\a,b} := \{\x \in \R^n : \a^\top \x = b\}, \ \ \ \H^{\leq}_{\a,b} := \{\x \in \R^n : \a^\top \x \leq b\}, \ \ \ \H^{\geq}_{\a,b} := \{\x \in \R^n : \a^\top \x \geq b\}.
$$
We say that $\H_{\a,b}$ (resp. $\H^{\leq}_{\a,b},\, \H^{\geq}_{\a,b}$) is a rational hyperplane (resp. rational halfspace) when $\a$ and $b$ are rational.

We use \(\langle \cdot \rangle\) to denote the binary encoding size of an object. For a rational number \(r = \frac{p}{q}\) written in lowest terms with integers \(p\) and \(q > 0\), the encoding size is defined as
\(
\langle r \rangle := 1 + \lceil \log_2(|p| + 1) \rceil + \lceil \log_2(q + 1) \rceil,
\)
which accounts for the binary representation of the numerator and denominator, as well as a sign bit.
For a vector \(\mathbf{v} = (v_1, \dots, v_n) \in \mathbb{Q}^n\), the encoding size is
\(
\langle \mathbf{v} \rangle := \sum_{i=1}^n \langle v_i \rangle.
\)
This convention ensures that the complexity of algorithms is measured in terms of the input bit-length, rather than its numeric magnitude.

For the convex sets we focus on, we will assume an exact separation oracle in the following definition.  
Note that we assume any output separating hyperplanes to be given with integral values, which could easily be obtained by scaling if instead they were encoded as rational values.  Such scaling had no significant influence on the encoding size of a hyperplane.

\begin{oracle}[(Exact) Separation Oracle for a Closed Convex Set]
\label{oracle:separationC}
A \emph{(exact) separation oracle} $\SEP{C}{D}$ for a closed convex set \( C \subseteq \mathbb{R}^n \) determines whether a given point from a domain $D$ is in $C$ or provides a hyperplane separating the point from $C$.

\end{oracle}

We will need to know the complexity of such oracles.  We define some general notation for this here.

\begin{definition}[Oracle Complexity Functions]
\label{def:oracle-complexity}
For any oracle $\mathcal{O}$, we associate two polynomials that are supposed take the encoding length of the input as an argument:
\begin{itemize}
    \item $\Phi^{\mathrm{out}}_{\mathcal{O}} : \mathbb{N} \to \mathbb{N}$, bounding the encoding length of the outputs of $\mathcal{O}$, and
    \item $\Phi^{\mathrm{time}}_{\mathcal{O}} : \mathbb{N} \to \mathbb{N}$, bounding the time (measured in the standard bit-complexity model, i.e., number of elementary steps) required to evaluate $\mathcal{O}$.
\end{itemize}
Specifically, for every query to $\mathcal{O}$ of encoding size at most $\tau$,
\[
\langle \mathcal{O}(\text{query}) \rangle \;\leq\; \Phi^{\mathrm{out}}_{\mathcal{O}}(\tau),
\]
and the evaluation runs in time at most $\Phi^{\mathrm{time}}_{\mathcal{O}}(\tau)$.
When the distinction is irrelevant, we write $\Phi_{\mathcal{O}}$ to denote a single polynomial that simultaneously bounds both.
\end{definition}

\begin{remark}
\label{rem:oracle-assumptions}
All oracles considered  are assumed to be 
\emph{deterministic}, and their output complexity functions 
$\Phi^{\mathrm{out}}_{\mathcal{O}}(\tau)$ are bounded by a polynomial in the input encoding size $\tau$.
\end{remark}

We say that an algorithm runs in \emph{oracle-polynomial time} if, assuming polynomial-time implementations of the oracles used, the algorithm runs in polynomial time.

The oracle above is efficiently computable in polynomial time for many well-known convex sets, even when the dimension is not fixed. 
For convex semi-algebraic sets, this can be done in polynomial time in fixed-dimension~\cite{khackiyan_integer_2000}.  See~\cite{Grtschel1993} for a thorough account of other oracles and how they relate to each other.
Separation oracles can even transition to intersections.

\begin{lemma}[Separation Oracles for Intersections]
\label{lemma:separation-intersection}
\label{rem:oracle-C-intersect-P}
Suppose \( C_1, C_2 \subseteq \mathbb{R}^n \) are closed convex sets presented by separation oracles \( \SEP{C_1}{D} \) and \( \SEP{C_2}{D} \), respectively, for some $D \subseteq \R^n$. 
Then we can construct a separation oracle \( \SEP{C_1 \cap C_2}{D} \) for \( C_1 \cap C_2 \).
\end{lemma}

\begin{proof}
Given a query point \( \x \in D \), we first invoke \( \SEP{C_1}{D} \) and \( \SEP{C_2}{D} \) on \( \x \). 
If both oracles confirm that \( \x \in C_1 \) and \( \x \in C_2 \), respectively, then return that \( \x \in C_1 \cap C_2 \).
 Otherwise, at least one of the oracles returns a valid inequality separating \( \x \) from its respective set. In that case, return that inequality as a valid separating hyperplane for \( C_1 \cap C_2 \).
\end{proof}

\subsection{Feasibility Problem}
We now formally define the core task of interest in an abstract setting.

 \begin{definition}[Feasibility Problem] \label{defn:feas-prob} Given a set \( S \subseteq \mathbb{R}^n \), the \emph{feasibility problem over \( S \)} is to either return a point \( \x \in S \), or assert that \( S = \emptyset \). \end{definition}

Depending on the structure of $S$, the corresponding feasibility problem can be solved efficiently.
This encourages us to define an oracle for this task, rendering it a tool for more advanced cases. 

\begin{oracle}[Feasibility Oracle]
\label{oracle:general-feas}
A \emph{feasibility oracle} $\FEAS{\mathcal Q}{D}$
for a domain $D \subseteq \R^n$ and a collection $\mathcal Q$ of sets $Q \subseteq \R^n$ takes on input a description of $Q \in \mathcal Q$ and determines if a point $\x \in Q \cap D$ exists, and if so, it returns $\x$.  
That is, $\FEAS{\mathcal Q}{D}$ solves the feasibility problem of over $Q \cap D$ for every $Q \in \mathcal Q$.
\end{oracle}

Note that opposed to the separation oracles being connected to one specific set, an oracle to solve the feasibility problem is supposed to be called for an entire class of sets, here denoted as $\mathcal Q$, 
since the process behind is typically available in a certain generality.
In the following, we will explore feasibility for both continuous and integer cases, i.e., $D = \Q^n$ or $\Z^n$, starting with $\mathcal Q$ consisting of convex sets described by separation oracles.

\subsubsection{Continuous Feasibility Testing via Inradius
}

While convex optimization is typically tractable, determining non-emptiness of a convex set using only a separation oracle is nontrivial. A standard approach is to either show that the set is non-empty, or it must not be full-dimensional because we can argue the volume is too small.  Hence, we can claim this convex set is contained in a hyperplane and then analyze it in a lower dimension.  For our purposes, we will need to know when two non-empty, open (and hence full-dimensional) convex sets intersect.  If they do, the intersection must have a positive volume.

The following result, adapted from~\cite[Theorem 3.3.3]{Grtschel1993}, provides a method for testing full-dimensionality using a separation oracle alone.

\begin{theorem}[Full-Dim Test with Inscribed Ellipsoid,~\cite{Grtschel1993}]
\label{thm:full-dim}
Let \( \epsilon, R > 0 \). There exists an oracle-polynomial time algorithm in \( \langle \epsilon \rangle, \langle R \rangle \) that takes as input a separation oracle $\SEP{C}{\Q^n}$ for a convex set \( C \subseteq B(\mathbf{0}, R) \) and returns one of the following:
\begin{enumerate}
    \item A full-dimensional ellipsoid \( E \subseteq C \), implying \( \intr(C) \neq \emptyset \).
    \item A declaration that \( \vol(C) \leq \epsilon \), implying \( C \) is not full-dimensional if \( \epsilon \) is chosen appropriately.
\end{enumerate}
\end{theorem}

The \emph{inradius} \( r \) of a convex set \( C \) is defined as the largest radius such that \( B(\x_0, r) \subseteq C \) for some $\x_0 \in C$. 
We assume a corresponding function \( r(C) \) exists such that, if \( C \) is full-dimensional, then \( B(\x_0, r(C)) \subseteq C \).

This allows setting \( \epsilon < \vol(B(\ve 0, r(C))) \) in the above theorem to test full-dimensionality.
We note the geometric inequality
\[
r \geq \frac{\vol(C)}{\vol(B(\ve 0, 1)) R^{n-1}},
\]
which provides a lower bound on the inradius in terms of the volume and outer radius of~\( C \). 
Though often stated informally, this result appears in \cite[Lemma 11]{Gribling2023} and underlies many bounds for convex bodies, particularly polytopes and semi-algebraic sets.
Several explicit inradius bounds are available:
\begin{itemize}
    \item For polytopes defined by integer matrices, \cite[Corollary 12]{Gribling2023} and \cite{Grotschel1988} provide volume-based bounds.
    \item For convex semi-algebraic sets defined by polynomial inequalities of degree \( d \), \cite[Proposition 2.5]{khackiyan_integer_2000} shows the inradius has encoding size \( \langle r \rangle = L d^{O(n)} \), where \( L \) is the encoding size of the input polynomials.
\end{itemize}

\begin{assumption}[Intersection Inradius]
\label{assumption:intersection-radius}
We assume that the inradius function behaves well under intersection, i.e.,
\[
\langle r(C_1 \cap C_2) \rangle \leq \operatorname{poly}(\langle r(C_1) \rangle, \langle r(C_2) \rangle ).
\]
That is, if \( C_1 \cap C_2 \neq \emptyset \), then the intersection retains a polynomially-bounded inradius.
\end{assumption}

This is critical for both pairwise intersections and intersections with halfspaces during cutting-plane procedures. For example, \cite[Lemma 3.1.35]{Grtschel1993} gives a lower bound on the volume of a full-dimensional polytope defined by integer data, which translates (via volume of balls) into an inradius bound.

\begin{corollary}
\label{cor:intersection-full}

Suppose \( C_1,\, C_2 \subseteq \R^n \) are closed convex sets presented by separation oracles $\SEP{C_1}{\Q^n}$ and $\SEP{C_2}{\Q^n}$ that satisfy \Cref{assumption:intersection-radius} and are contained in a ball of radius \( R > 0\). 
Then in oracle-polynomial time, we can test whether \( C_1 \cap C_2 \) is full-dimensional. 
Moreover, the same algorithm determines whether \( \intr(C_1) \cap \intr(C_2) = \intr(C_1 \cap C_2) \) is nonempty.

\end{corollary}

\begin{proof}
We first construct a separation oracle for \( C_1 \cap C_2 \) using \Cref{lemma:separation-intersection}. 
Then we apply \Cref{thm:full-dim} with \( \epsilon < \vol(B(\ve 0, r)) \), where \( r = r(C_1 \cap C_2) \). 
By \Cref{assumption:intersection-radius}, the encoding size of \( r \) remains polynomial, so the test runs in oracle-polynomial time.
\end{proof}

\begin{remark}
In special cases, such as intersections of balls, full-dimensionality can be verified directly; see \Cref{remark:sphere-intersection}.
\end{remark}

\subsubsection{Integer Feasibility Oracles}
\label{subsec:feasibility-oracles}

For many types of collections of bounded convex sets $\mathcal C$, it is well known that an oracle $\FEAS{\mathcal C}{\Z^n}$ can be implemented in fixed-dimension $n$, 
typically with a fixed parameter tractable result (where $n$ does not appear as an input of the encoding of the sets in the complexity).
This was initially done for polyhedra by Lenstra~\cite{lenstra_integer_1983}, but has been extended to other contexts. 
 We refer to~\cite{dadush2012integer} for a thorough treatment using separation oracles $\SEP{C}{\Q^n}$. 
 See  \cite[Theorem 7.1.1]{dadush2012integer} and also \cite{rothvoss-subspace} for improved complexity results.
Such an algorithm for integer feasibility on a bounded region is easily transferred to an algorithm for optimization over integer points via a bisection-based search over objective values.

 As pointed out in \cite{oertel}, $\FEAS{\mathcal C}{\Z^n}$ can be solved using oracle queries only on integer points, i.e., $\SEP{C}{\Z^n}$.  
 Although,  as far as we know, this currently requires a different approach and results in a worse complexity; though it is still polynomial time in fixed dimension $n$.
 We state an adapted version of this result and for completeness, we provide a proof of this results in \Cref{sec:complexity-integer-queries} that mimics a proof from Basu~\cite{Basu2022}.
 Querying the separation oracle only on integer points will be crucial for our results.

\begin{theorem}[Convex Integer Feasibility and Optimization~\cite{oertel}, \cite{Basu2022}]
\label{thm:convex-integer} Let $\mathcal C$ be the  collection of convex sets $C \subseteq B_\infty(\mathbf{0},R)$ presented by a separation oracle $\SEP{C}{\Z^n}$.  
Then there is an algorithm $\FEAS{\mathcal C}{\Z^n}$ 
that runs in time 
\begin{equation*} 
2^{O(n\log(n))}\operatorname{poly}(\langle R \rangle, \Phi_{\SEP{C}{\Z^n}}(n\langle R \rangle ) )
\end{equation*}
for each $C \in  \mathcal C$.

If feasible, an optimal solution to $\max/\min \{\c^\top \x : \x \in C \cap \Z^n\}$ can be computed in the same order of time and space complexity.

\end{theorem}

Let $\PC$ denote the collection of reverse convex sets given as $P \setminus C$ where $P,\, C\subseteq \R^n$, $P$ is a polytope given by an inequality description, and $C$ is a convex set.  
As demonstrated in~\cite{cubic-integer-plane} the oracle $\FEAS{\PC}{\Z^n}$ can be implemented in polynomial time in fixed-dimension.  
The key is to compute the integer hull $P_I$ of $P$ in polynomial time in fixed-dimension using \cite{cook_integer_1992,hartmann-1989-thesis}.  
The integer hull of a rational polyhedron has at most $2m^n(6n^2 \phi)^{n-1}$ many vertices where $m$ is the number of inequalities and $\phi$ is the binary encoding size of the polyhedron.  
This approach is based on the following property.
\begin{quote}
\emph{Since $C$ is convex, the set $(P \setminus C) \cap \Z^n$ is non-empty if and only if there exists a vertex of the integer hull of $P$ that is not in $C$. } 
\end{quote}
Thus, we can test for feasibility in this setting by enumerating the vertices of the convex hull $\conv(P \cap \Z^n)$.  
Notably, this algorithm only requires membership access to $C$.  
Unfortunately, this property does not carry over to removing multiple convex sets.
\begin{example}
Consider the following example.\\
\begin{minipage}{0.3\textwidth}
\begin{center}
   \begin{tikzpicture}

\fill[yellow, opacity=1, draw = black] (1,0) -- (0,1) -- (1,2) -- (2,2) -- (2,1) -- cycle;
\foreach \x in {0,1,2}
    \foreach \y in {0,1,2}
        \node at (\x,\y) [circle,fill,inner sep=1pt] {};

\fill[blue, opacity=0.3] (1,1.7) ellipse [x radius=1.5, y radius=0.5, rotate=35];

\fill[green, opacity=0.3] (1.7,1) ellipse [x radius=1.5, y radius=0.5, rotate=65];

\end{tikzpicture}
\end{center}
\end{minipage}
\begin{minipage}{0.7\textwidth}
Let $P = \conv(\{(1,0), (0,1), (1,2), (2,2), (2,1)\})$, $C_1 \supseteq \conv(\{ (0,1), (1,2) , (2,2) \})$ and $C_2 \supseteq \conv(\{(1,0), (2,1), (2,2)\})$.
Then, as demonstrated,  $C_1 \cup C_2$ contains the vertices of the integer hull of $P$ (which is the same as vertices of $P$ in this case).  However, $P\setminus (C_1 \cup C_2)$ contains the integer point $(1,1)$.
\end{minipage}
\end{example}

Our goal is to utilize oracles $\FEAS{\mathcal C}{\Z^n}$ and $\FEAS{\PC}{\Z^n}$ for solving non-convex integer optimization problems.  
We focus on results in general (fixed) dimensions $n$.
To this effect, we generalize a notion from~\cite{cubic-integer-plane}.

\begin{definition}
\label{def:division-description}
Given a set $S$ and a box $B= B_\infty(\ve 0, R)$, a \emph{\CCPsub} of $S$ on $B$ is a list of rational polyhedra $P_i$, $i = 1, \dots, l_1$, $Q_j$, $j = 1, \dots, \ell_2$, such that
\begin{enumerate}[(i)]
\item $P_i \cap S$ is convex for $i=1, \dots, \ell_1$,
\item $Q_j \setminus S$ is convex for $j=1, \dots, \ell_2$, and
\item 
\begin{equation*}
B\cap \Z^n = \left(\bigcup_{i=1}^{\ell_1} P_i \cup \bigcup_{j=1}^{\ell_2} Q_j
\right) \cap \Z^n .
\end{equation*}
\end{enumerate}
\end{definition}

\begin{theorem}[Integer feasibility given a subdivision]
\label{thm:dd-min}
Suppose $S \subseteq B_\infty(\ve 0, R)$ is presented with a {\CCPsub} on this box.
Further, suppose we are provided with oracles $\FEAS{\mathcal C}{\Z^n}$ and $\FEAS{\PC}{\Z^n}$. 
Then, if the dimension $n$ is fixed, we can solve the feasibility problem over $S\cap \Z^n$ in oracle-polynomial time .
\end{theorem}

\begin{proof}
For each $i=1, \dots, \ell_1$, apply $\FEAS{\mathcal C}{\Z^n}$ to $P_i \cap S$.
Next, for each $i=1, \dots, \ell_2$, apply $\FEAS{\PC}{\Z^n}$ with  $P = Q_j$ and $C = Q_j \setminus S$ which solves the integer feasibility problem over $P \setminus C= Q_j \setminus (Q_j \setminus S) = Q_j \cap S$. 

If any oracle call returns a feasible point, then we answer YES to the feasibility problem and return that point.
Otherwise, we answer NO to the feasibility problem.

As discussed, the oracles used can be implemented in polynomial time in fixed-dimension.  
Thus we obtain a complexity that is polynomial-time in fixed-dimension (and polynomial in $\ell_1, \ell_2, d$ and the maximal encoding size $\phi$ of the polyhedra $P_i$, $Q_j$).
\end{proof}

\subsection{Tools}

We now revisit how the separation–optimization equivalence can be used to construct oracles for convex operations, assuming inradius bounds as described above.

\begin{lemma}[Separation Oracle for a Convex Union]
\label{rem:oracle-union}
Let \( C_i \subseteq \mathbb{R}^n \) be closed convex sets for \( i = 1, \dots, m \), and suppose that the union \( C := \bigcup_{i=1}^m C_i \) is convex. If separation oracles \( \SEP{C_i}{\Q^n} \) are available for each \( C_i \), then a separation oracle \( \SEP{C}{\Q^n} \) can be realized using these oracles.
\end{lemma}

\begin{proof}
Since \( C \) is convex and expressed as a finite union of convex sets, it must in fact coincide with the convex hull of their union: \( C = \operatorname{conv}\left( \bigcup_{i=1}^m C_i \right) \). 
This structure allows us to construct a linear optimization oracle over \( C \) by optimizing over each \( C_i \) independently and returning the point achieving the best objective value, i.e.,
\[
\max_{\x \in C} \mathbf{c}^\top \x = \max_{i=1, \dots, m} \left\{ \max_{\x \in C_i} \mathbf{c}^\top \x \right\}.
\]
By the separation–optimization equivalence (e.g., \cite{Grotschel1988}), and assuming bounds on inradius and circumscribing radius, this suffices to define a separation oracle for \( C \) in oracle-polynomial time.
\end{proof}

\begin{remark}
This result assumes access to linear optimization over each \( C_i \), which may itself rely on the equivalence between optimization and separation for each \( C_i \). 
The inradius of \( C = \bigcup_i C_i \), if convex, can be lower bounded by the minimum inradius among the~\( C_i \), provided the union is structured.
\end{remark}

The result is due to a variant of the ellipsoid method that achieves \emph{ellipsoid rounding}, i.e., finds an ellipsoid inscribed in $C$ such that a scaling of the ellipsoid circumscribes $C$.  Similar results can be obtained via polyhedral approaches such as in~\cite{Anstreicher1999}.

In the remainder of this subsection, we introduce useful lemmas that are leveraged throughout the article. 
To keep the presentation concise, we omit the routine proofs of the first two statements.

\begin{lemma}
\label{obs:intr-intersection}
    Let $C_1, C_2\subseteq \R^n$.  Then
    \(
\intr(C_1) \cap \intr(C_2) = \intr(C_1 \cap C_2).
    \)
\end{lemma}

\begin{lemma}
\label{lem:interior-intersection-full-dim}
Let \( C_1, C_2 \subseteq \mathbb{R}^n \) and suppose \( C_2 \) is full-dimensional. Then
\[
\intr(C_1) \cap C_2 \neq \emptyset \quad \Longleftrightarrow \quad \intr(C_1 \cap C_2) \neq \emptyset.
\]
\end{lemma}

The following lemma allows us to relax a polytope to one containing the same integer points, but all integer points are now are the interior of the relaxed polytope.

\begin{lemma}[Integer Point Preserving Relaxation]
\label{lem:interior-intersects-lattice}
Let 
\(
P = \{\, \x \in \R^n : A\x \leq \b \,\},
\)
where \(A \in \Z^{m \times n}\) and \(\b \in \Z^m\), be a full-dimensional rational polyhedron.  
Define
\[
P' := \{\, \x \in \R^n : A\x \leq \b + \tfrac{1}{2} \mathbf{1}\,\}. 
\]
where $\mathbf{1}$ is the vector of all 1's.
Then \(P\) and \(P'\) contain the same integer points, with \(P \subseteq \intr(P') \subseteq P'\).
Moreover, for any full-dimensional, closed, convex set \( C \subseteq \R^n \),
\[
\intr(C) \cap P \cap \Z^n = \intr(C \cap P') \cap \Z^n.
\]
Finally, if \( \intr(C) \cap P \cap \Z^n \neq \emptyset \), then \( \intr(C) \cap P' \) is full-dimensional.
\end{lemma}

\begin{proof}
 Since each inequality \(\a^\top \x \le b\) of \(P\) has integer right-hand side,
every integer point \(\z \in \Z^n\) satisfies either \(\a^\top \z \le b\) or 
\(\a^\top \z \ge b+1\). Thus
\[
\a^\top \z \le b \;\; \Longleftrightarrow \;\; \a^\top \z \le b + \tfrac{1}{2}.
\]
Hence \(P \cap \Z^n = P' \cap \Z^n\). 
It is also clear that \(P \subseteq \intr(P')\) and that
\(P'\) is rational.

Next, we show the equivalence of sets.

($\subseteq$) Let \(\z \in \intr(C) \cap P \cap \Z^n\). Then \(\z \in \intr(C)\) and, since \(\z \in P \subseteq \intr(P')\), Lemma~\ref{obs:intr-intersection} gives
\[
\z \in \intr(C) \cap \intr(P') = \intr(C \cap P').
\]
Thus, \(\z \in \intr(C \cap P') \cap \Z^n\).

($\supseteq$) Conversely, let \(\z \in \intr(C \cap P') \cap \Z^n\). 
Then \(\z \in \intr(C) \cap \intr(P')\) by~\Cref{obs:intr-intersection}.
As shown above, every integer point of \(P'\) is also in \(P\), hence \(\z \in P\). 
Therefore \(\z \in \intr(C) \cap P \cap \Z^n\).

\medskip
Lastly,  if \(\intr(C) \cap P \cap \Z^n \neq \emptyset\), then by the above, we have \(\intr(C \cap P') \cap \Z^n \neq \emptyset\). 
Since both \(C\) and \(P'\) are full-dimensional closed convex sets and the intersection of their interiors is nonempty, 
the intersection is also full-dimensional, so \(\intr(C \cap P')\) is nonempty and full-dimensional. 
Consequently, \(\intr(C) \cap P'\), which contains \(\intr(C \cap P')\), is also full-dimensional.
\end{proof}
The role of the relaxed polyhedron \(P'\) is to ``round out'' the facets of \(P\), i.e.,
every facet inequality is shifted outward by \(\tfrac{1}{2}\). This guarantees that 
\(P \subseteq \intr(P')\) while preserving all integer points 
(\(P \cap \Z^n = P' \cap \Z^n\)). 
As a result, intersections of the form 
\(\intr(C) \cap P \cap \Z^n\) can be expressed equivalently as 
\(\intr(C \cap P') \cap \Z^n\), which is often more convenient to analyze.
In particular, it allows us to apply interior arguments to the intersection
\(C \cap P'\) without changing the underlying set of integer feasible points.

\subsection{Separating Convex Sets}
A key tool that we will use is to separate two convex sets that don't intersect full-dimensionally.  
Ideally, we would like an algorithm that, given closed convex sets  $C_1$, $C_2$, detects if $\intr(C_1 \cap C_2) = \emptyset$,
and if so, returns a hyperplane $H_{\a,b}$ that separates $C_1$  and~$C_2$. 

For two convex sets $C_1, C_2$, we say that they are \emph{weakly separable} if there exists $(\a,b)$ such that $C_1 \subseteq H^{\leq}_{\a,b}$, $C_2 \subseteq H^{\geq}_{\a,b}$. 
We formalize this as an oracle.

\begin{oracle}[Weak Convex Set Separation]
 \label{oracle:separation}
 For a collection of  sets $\mathcal Q$, the oracle $\WEAK{\mathcal Q}{D}$ takes as input two sets $Q_1, Q_2 \in \mathcal Q$ and  either outputs a hyperplane $H_{\a,b}$ that weakly separates $Q_1 \cap D$ and $Q_2 \cap D$ or confirms that no such weakly separating hyperplane exists.
 
 \end{oracle}
 
 However,  for weakly separable convex sets $C_1,\, C_2$ that are presented only via exact separation oracles, it can be  difficult to exactly confirm separation of the two convex sets.  
 In~\cite[Theorem 4.7.1]{Grtschel1993}, they show how this separation can be done via oracles, but only provide an $\varepsilon$-approximate separation guarantee. In later sections, we will establish this separation over just the integer points (\Cref{lem:separation-integer-hulls}).

\subsection{Reverse Convex Sets}
An important result that allows for understanding reverse convex sets (over continuous variables) is that the convex hull is a polytope.  We will use this fact in our proofs.

\begin{theorem}[\cite{Hillestad1980}]
\label{thm:reverse_convex_polytope}
Let $P \subseteq \mathbb{R}^n$ be a polytope and let $C_1, \ldots, C_m$ be closed convex sets. Then
$$
Q = \operatorname{conv}\left(P \backslash\left(\bigcup_{i=1}^m \operatorname{int}\left(C_i\right)\right)\right)
$$
is a polytope.
\end{theorem}

\begin{remark}
As mentioned in~\cite{Jacobsen2008}, a formulation with multiple reverse convex constraints can be lifted by one variable to a problem with a single reverse convex constraint.  However, this also introduces an additional convex constraint. 
In our setting, we want to ensure that the convex constraints are polyhedral so that we can efficiently study the integer hull of these points.  Thus, we cannot apply this great trick of rewriting to a single reverse convex constraint.
\end{remark}

\subsection{Hyperplane Arrangements}

A \emph{hyperplane arrangement} $\mathcal{H}$ is a finite collection of hyperplanes $\{\H_{\a_i,b_i}\}_{i=1}^d$ in $\R^n$.
The \emph{cells} $\mathcal P$ of a hyperplane arrangement $\mathcal{H}$ are the set of polyhedra defined as
$
P = \cap_{i=1}^d \H^{*_i}_{\a_i, b_i}
$
for $*_i \in \{\leq, \geq\}$. The set of maximal cells, $\mathcal P_{\max} \subseteq \mathcal P$ is the subset of cells that are full-dimensional.  
Alternatively, these can be viewed as the closure of connected components of the complement of the union of the hyperplanes in $\mathcal{H}$, i.e.,
$
\mathbb{R}^n \setminus \bigcup_{i=1}^d H_{\a_i,b_i}.
$

The following result is well-known and we do not know the exact origin.  See, for example, \cite[Section 1.2]{Edelsbrunner1987}.
\begin{lemma}[Number of cells in hyperplane arrangement]
\label{lem:hyperplane-arrangement-bound}
Let $\mathcal{H}$ be a hyperplane arrangement of $d$ hyperplanes in $\R^n$.  Let $\mathcal P_{\max}$ be the maximal cells of $\mathcal H$.  Then $|\mathcal P_{\max}| = O(d^n)$.
\end{lemma}

Various algorithms can enumerate the maximal cells of hyperplane arrangements.  We cite the following theorem for an output sensitive approach.  
\begin{lemma}[{\cite[Theorem 1]{sleumer_output_1998}}]
\label{lem:hyperplane-arrangement}
Let $\mathcal H$ be a hyperplane arrangement in $\R^n$ with $d$ hyperplanes.  The maximal cells $\mathcal P_{\max}$ can be enumerated in time 
$$
O\left(\left(n \cdot \operatorname{lp}(d, n)+d \cdot n^2+d \cdot \operatorname{lp}(l, n)\right) \cdot|\mathcal P_{\max}|\right),
$$
where $l$ is the maximum number of neighbors of a cell and $\operatorname{lp}(d, n)$ indicates the time (or space) complexity needed for solving a $d \times n$ linear program.   If $\mathcal H$ is rational with binary encoding size $\phi$, then this time complexity is bounded by
$$\mathrm{poly}(n,d,\phi)O(d^n).$$
\end{lemma}

\section{Hard and Simple Cases}
\label{sec:hard-simple-cases}

We seek complexity results for \eqref{eq:main}  for sets where $K_i$ and $C_i$ are convex sets.  We will use the technique of \CCPsubs as the key approach.
We use $K$ and $C$ when talking about convex sets and use $P$ and $Q$ instead when these are polyhedra.  
Certain versions of this may be difficult.

\subsection{Hard Cases: Unspecific Bounds}
We give two examples of hard cases;  both are derived from examples in~\cite{cubic-integer-plane} in a different format to show difficulties in the context of minimizing degree 4 polynomials over the integer points in a polytope. \\

\begin{minipage}{0.35\textwidth}
\hspace{-1cm}
\input{AN1.tex}
  \end{minipage}
\begin{minipage}{0.6\textwidth}
\begin{example}[NP Hard $(K \setminus C) \cap \Z^2$]
It is NP-Hard to determine if $K\setminus C$ contains an integer point when $K$ and $C$ are semialgebraic sets described by even one inequality of degree 2 in dimension 2 along with some affine inequalities.  This reduces to the NP-Hard Problem AN1, which essentially relies on detecting integer points on a quadratic curve~\cite{manders_np_1976}.\\
In particular, \cite{loera_integer_2006} shows that problem AN1 can be phrased as minimizing $(ax - b - y^2)^2$ over the integer points in a box.
If the optimal solution is 0, then there is a \texttt{YES} answer to the AN1 problem. Otherwise, there is a \texttt{NO} answer.  
This is then equivalent to testing the integer feasibility of $K \setminus C$ where 
$$
K = \{(x,y) \in B : ax - b - y^2 \geq 0\},
$$
$$
C = \{(x,y) \in B : ax - b - y^2 \geq  1\},
$$
where $B = [\ell_1, u_1] \times [\ell_2, u_2]$ for some values $\ell_i, u_i \in \R$.  In particular, $K$ and $C$ are bounded.  Then, given the correct choices of bounds, we can assert a \texttt{YES} answer to AN1 if the set is feasible over the integers, and \texttt{NO} otherwise.
\end{example}
\end{minipage}

\vspace{0.3cm}
\hspace{-0.6cm}\begin{minipage}{0.6\textwidth}
\begin{example}[Large Solutions in Unbounded Case]
\label{ex:unbounded}
The Pell equation is 
$(x^2 - N y^2)^2 = 1$ and is to be solved for $x,y \in \Z_{\geq 1}$.  
When $N$ is a large integer $N = 5^b$ with $b \in \mathbb{N}$, it is known that solutions $x,y$ can be extremely large, i.e., $|x|,|y| =\Omega( 2^{2^b})$~\cite{lagarias_computational_1980}. 
Thus, writing down solutions can be difficult.

This set can be realized as $P \setminus (C_1 \cup C_2)$ where \\
$P = \{(x,y) : x,y \geq 1\}$,\\
$C_1 = \{(x,y) : x^2 - N y^2 \geq 2\}$ and \\$C_2 = \{(x,y) : x^2 - Ny^2 \leq -2\}$.
\end{example}
\end{minipage}
\begin{minipage}{0.35\textwidth}
\hspace{.1cm}
\input{pell-pic.tex}
  \end{minipage}

\subsection{Simple Case: All Removed Sets are Polyhedra }

If all sets being removed are polyhedra, we can then identify a decomposition by looking at the appropriate hyperplane arrangement.
\begin{theorem}
\label{thm:removing-polytopes}
    Let $K$ be a convex set and let \(Q_1, \dots, Q_m\) be a set of rational polyhedra in $\R^n$ defined by inequality representations.  Then there exist polyhedra $T_1, \dots, T_t$ such that 
    \begin{equation}
        \label{eq:k-minusQs}
    \left(K \setminus \left(Q_1 \cup \dots \cup Q_m\right)\right) \cap \Z^n = \left(\bigcup_{s=1}^t (K \cap T_s)\right) \cap \Z^n.
    \end{equation}
Here, $t = O(d^n)$, where $d$ is the total number of inequalities needed to describe the polyhedra $Q_1, \dots, Q_m$ and the encoding sizes of $T_s$ are bounded by a polynomial in the encoding sizes of $Q_1, \dots, Q_m$.
\end{theorem}
\begin{proof}
Consider the case $Q_1 \cup \dots \cup Q_m = \R^n$. 
Since $K \setminus (Q_1 \cup \dots \cup Q_m) = \emptyset$, it suffices to set $t = 1$ and define $T_1$ as some polyhedron without any integer points in its interior, for instance, the empty set.

Turning to the case $Q_1 \cup \dots \cup Q_m \subsetneq \R^n$, 
suppose there is a rational polyhedron $Q = \{\x : A\x \leq \b\}$ with $A$ and $\b$ rational.
Without loss of generality, $A$ and $\b$ can even be considered integral. 
By~\Cref{lem:interior-intersects-lattice},  $Q' = \{\x : A\x \leq \b + \tfrac{1}{2} \mathbf{1}\}$ satisfies
$$
Q \cap \Z^n = Q' \cap \Z^n.
$$
Furthermore, since $A\x$ is an integer vector for any $\x \in Q \cap \Z^n$, we see that $\x \in \intr(Q')$ since it satisfies every constraint with strict inequality.
Hence, considering the rational polyhedra $Q_i$ in the statement of this theorem, without loss of generality, we will assume that $Q_i \cap \Z^n \subseteq \intr(Q_i)$ by reassigning $Q_i \leftarrow Q'_i$.

Now, construct a hyperplane arrangement $\mathcal{H}$ comprising all the inequalities that define  $Q_1, \dots, Q_m$. Using Lemma~\ref{lem:hyperplane-arrangement}, we enumerate the maximal cells of the hyperplane arrangement $\mathcal{H}$, which can be done in a number of steps that is polynomial in the size of the input. Next, let $\{T_s\}_{s = 1}^t$ be the subset of maximum cells such that $\intr(T_s) \cap (Q_1 \cup \cdots \cup Q_m) = \emptyset$.  Since the maximum cells partition the integer points $\Z^n$, we have that 
$$
\left ( \bigcup_{s=1}^t T_s \right) \cap \Z^n = \Z^n \setminus \left(Q_1 \cup \cdots \cup Q_m \right).
$$
Thus, \eqref{eq:k-minusQs} follows.  

\end{proof}

The above theorem produces a number of cells that are exponential in the dimension. 
Unfortunately, this is to be expected.  
Consider the example in the following remark.

\begin{remark}
    For $m = n$, it is unlikely to obtain a sub-exponential complexity in the dimension even if we relax the integrality constraints.  In particular, even with choosing $Q_i = \{\x :  \varepsilon \leq x_i \leq 1-\varepsilon\}$ for some $\varepsilon > 0$ and suppose $K \subseteq [0,1]^{n}$.  Then we see that $\lim_{\varepsilon \to 0} [0,1]^n \setminus (Q_1 \cup \dots \cup Q_m) = \{0,1\}^n$.  Thus, choosing $\varepsilon$ sufficiently small recovers approximately $K \setminus (Q_1 \cup \dots \cup Q_m) \approx K \cap \{0,1\}^n$.  That is, we can approximate any binary integer program, and thus we expect that runtime for even detecting a feasible point to be $\Omega(2^n)$.
\end{remark}

\begin{corollary}
\label{cor:Ci-polyhedral}
The feasibility problem over~\eqref{eq:main} can be solved in polynomial time in fixed-dimension when all sets $C_i$ are polyhedral.
\end{corollary}
\begin{proof}
Let $C_i$ denoted by $Q_i$ as they are  polyhedral.  
We test integer feasibility on each $S_j := K_j \setminus (Q_1 \cup \dots \cup Q_m)$, $j = 1,\dots,\ell$, by means of a \CCPsub\ of $S_j$ on the box containing $K_j$.
In particular, we use~\Cref{thm:removing-polytopes} to compute polyhedra $T_s$, which fulfill 
$\intr(T_s) \cap (Q_1 \cup \cdots \cup Q_m) = \emptyset$.
Thus, $S_j  \cap T_s = K_j \cap T_s$ is convex and the $T_s$ serve as the ``$P_i$" in~\Cref{def:division-description}.
The ``$Q_j$" from the same definition are exactly the ones subtracted from $K_j$.

The \CCPsub\ allows to apply Theorem~\ref{thm:dd-min} to detect integer feasibility.  
Practically, only the oracle for convex integer feasibility (\Cref{thm:convex-integer}) is needed, since only the $T_s \cap K_j$ have to be checked.
If the integer feasibility test on one $S_j$ asserts a point $\x^*$, we return it.  Otherwise, we assert that there is no such solution.
\end{proof}

\section{Polytope Minus Convex Sets}
\label{sec:main-continuous}
We now focus on the setting 
\begin{equation*}
    P \setminus (C_1 \cup C_2\cup \dots \cup C_m),
\end{equation*} 
where $P$ is a polyhedron while the $C_i$ may be more general convex sets.  
As shown in \Cref{ex:unbounded}, if $P$ and the convex sets $C_i$ are unbounded, we may have complexity issues in even writing down solutions. 
If we assume that the integer hull of $P$ is unbounded, but $C_i$ are bounded, then trivially there is a feasible integer point attained by finding a sequence of integer points in $P$ that diverge.
Thus, this case can be solved conveniently and will not be discussed further.

\textbf{Boundedness Assumption.} Thus, we will assume that  $P$ is a bounded polyhedron (i.e., polytope).  In this setting, we can assume also that $C_i$ are bounded as well by redefining them as $C_i \leftarrow P \cap C_i$.  

In our main result, we will assume additional properties on the sets $C_i$.
We begin with a remark about why a trivial approach might not make sense.  
\begin{remark}
Perhaps a simple approach to solve this kind of feasibility problem would be to simply convert $C_i$ to polyhedra.  Indeed, if we are only concerned with feasibility over the integer points, we could instead replace $C_i$ by $Q_i$ where
$$
Q_i = \conv(C_i \cap P \cap \Z^n).
$$
However, even in the case where $C_i$ is a Euclidean Ball, any facet or vertex description of $Q_i$ would be expected to be exponential in the size of the input; in particular, in the radius of the ball~\cite{Brny1998}.
\end{remark}

\subsection{\HC}
\label{subsec:BHC}
Since we may need to deal with complicated integer hulls, we now focus on convex sets with the following property.  We write the property here without requiring the sets to be convex, but we will always apply it to convex sets.
\begin{definition}[\HC]
\label{defn:bhc}
    Let $C_1, \dots, C_m \subseteq \R^n$.  We say that $C_1, \dots, C_m$ have a \emph{\HC} if there exists a finite collection $H^1, \dots, H^d$ of hyperplanes such that for every pair $(i,j)$ with $1 \leq i < j \leq m$, we have that ${\bd C_i \cap \bd C_j \subseteq \bigcup_{s = 1}^d H^s}$.
\end{definition}

We motivate this definition by an example.
 We say the Euclidean ball $B(\c, r)$ is rational if $\c,\,r$ are both rational and we refer to its boundary -- the sphere -- as $\partial B(\c,r)$.

\begin{proposition}
\label{lem:spheres}
    The intersection of two spheres, which are not identical, is contained in a hyperplane.  
    If the spheres are rational, then the hyperplane is also rational.
\end{proposition}
\begin{proof}
Without loss of generality, we assume that one sphere is the boundary of the unit ball, while the other sphere is given by a center $\c \in \R^n$ and a radius $r$.

Then the intersection 
\begin{align*}
\bd B(\ve 0, 1) \cap \bd B(\c, r) &= \{ \x \in \R^n : \|\x\|_2 = 1, \|\x-\c\|_2 = r\}\\
& = \{ \x \in \R^n : \sum_{i=1}^n x_i^2 = 1, \sum_{i=1}^n x_i^2 - 2\sum_{i=1}^n c_i x_i + \sum_{i=1}^n c_i^2 = r^2\}\\
& = \{ \x \in \R^n : \sum_{i=1}^n x_i^2 = 1, 1 - 2\sum_{i=1}^n c_i x_i + \sum_{i=1}^n c_i^2 = r^2\}.
\end{align*}
Thus, the desired hyperplane is given by $H = \{\x \in \R^n: 1 - 2\sum_{i=1}^n c_i x_i + \sum_{i=1}^n c_i^2 = r^2\}$, 
which is rational provided that $\c$ and $r$ are rational.
\end{proof}
Note that the set $\bd B(\ve 0, 1) \cap \bd B(\c, r) $ from \Cref{lem:spheres} might be lower-dimensional or even empty.  This can happen if $B(\c, r) \subseteq B(\ve 0, 1)$ or $B(\ve 0, 1) \cap B(\c, r)$ is empty or a single point.

\begin{corollary}
\label{cor:ball-bhc}
 Let $B(\c^i, r_i)$, $i=1, \dots, m$, be Euclidean balls.  Then they have a \HC  with ${m \choose 2}$ hyperplanes.
\end{corollary}
\begin{proof}
    Apply \Cref{lem:spheres} to every pair of balls.
\end{proof}

\begin{remark}[Sphere intersection]
\label{remark:sphere-intersection}
Note that it is not difficult to derive characterizations on when two spheres intersect.  
The spheres intersect if and only if
$ r - 1\leq  \|\c\|_2 \leq r +1$.
Clearly, when $\|\c\|_2 = r+1$,  the spheres touch in a single point, and when $\|\c\|_2 = r-1$,  the second ball contains the unit ball.  In all other feasible intersection cases, their intersection is a $(n-1)$-sphere  defined by slicing the unit sphere by a hyperplane.
\end{remark}
Balls even fall into a class showing a special type of \HC, namely the \emph{ideal} one, as described in the following.

\subsection{Ideal Covers and Complicating Examples}

Let $C_i$, $i =1,\dots,m$ be convex sets that satisfy \HC with a collection of hyperplanes $\mathcal H$.
Then, imagine the (maximal) cells resulting from $\mathcal{H}$.
A first thought is that there appears only one connected convex set (possibly resulting from unions) in each such cell.
Although this is not necessarily true in general, there are nice criteria that this may happen.
\begin{definition}[Ideal Covers]
\label{defn:ideal-boundary-cover}
    Let $C_1, C_2$ be full-dimensional closed convex sets.  Let $H$ be a hyperplane such that $\bd C_1 \cap \bd C_2 \subseteq H$.  We say that $H$ is an \emph{ideal boundary cover} of $C_1, C_2$ if $H = \aff(\bd C_1 \cap \bd C_2)$. 

\label{defn:ideal-BHC}
    Further, let $\{C_i\}_{i=1}^m$ be full-dimensional closed convex sets.  We say that $\mathcal H$ is an \emph{Ideal \HC} if it is a \HC such that  every pair of intersecting $C_i,\, C_j$ has an ideal boundary cover $H \in \mathcal H$. 
\end{definition}

\begin{figure}[h!]
\begin{center}
\subfigure[]{
\includegraphics[scale = 0.13]{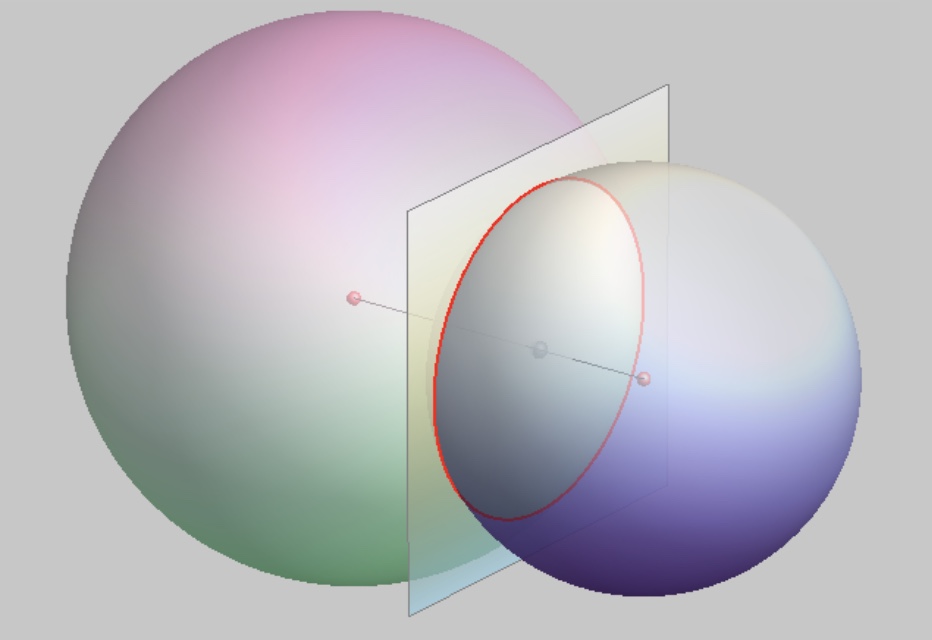} \hspace{1cm}
}
\subfigure[]{
\begin{tikzpicture}
    \fill[blue, opacity = 0.2] (-1,0) circle (1);

    \fill[green, opacity = 0.2] (1,0) circle (1);

    \draw[] (0,-1.5) -- (0,1.5);
    \draw[thick, dotted] (1.7,1) -- (-1.7,-1);
    \draw[dashed] (-0.5,1.5) -- (0.5,-1.5);
\end{tikzpicture}
}
\hspace{1cm}
\subfigure[]{
\begin{tikzpicture}[scale = 1.0]
    \fill[blue, opacity = 0.2] (0,0) -- (1,0) -- (1,1) arc[start angle=0, end angle=180, radius=0.5] -- (0,0);
    \fill[green, opacity=0.2] (0,0)--(0,1) -- (1,1) arc[start angle=90, end angle=-90, radius=0.5] -- (1,0) -- (0,0);

    \draw[dotted] (-0.5,0) -- (2.5,0); %
    \draw[dashed](0,-0.5) -- (0,2.5); %
    \draw(2.5,-0.5) -- (-0.5,2.5); %
\end{tikzpicture}
}
\end{center}
\caption{Types of Boundary Hyperplane Covers.  (a) Ideal type of boundary hyperplane cover since it results in exactly one convex set of relevance on each side of the hyperplane.  (b) The hyperplane cover might not be unique: we display three different hyperplanes that satisfy the definition.  In this case, the sets are weakly separable, but that might not be obvious from the start.  (c) Complicated boundary hyperplane cover with three hyperplanes.  These hyperplanes cover the intersections of the convex sets.  Some cells of this hyperplane arrangement do not intersect the convex sets.  The interior of the formed triangle intersects both of them.  In this case, the resulting intersection is a union of convex sets that is a convex set.  Lastly, the upper right cell has both convex sets, but they are weakly separable on this cell.}
\label{fig:types-bhc}
\end{figure}
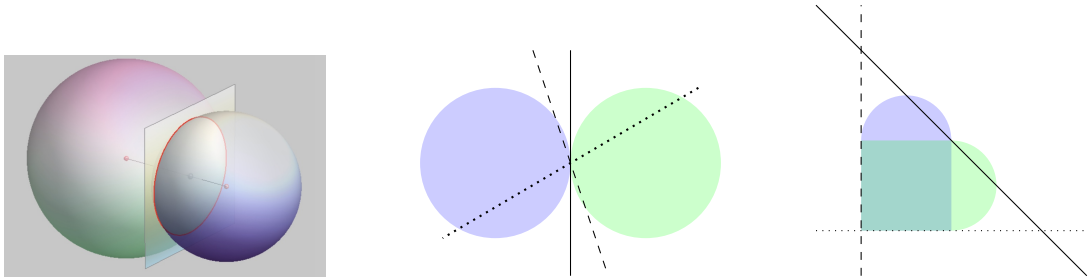

\begin{figure}[h!]
    \centering
\includegraphics[scale = 0.1, trim = 700pt 0pt 0pt 0pt, clip]{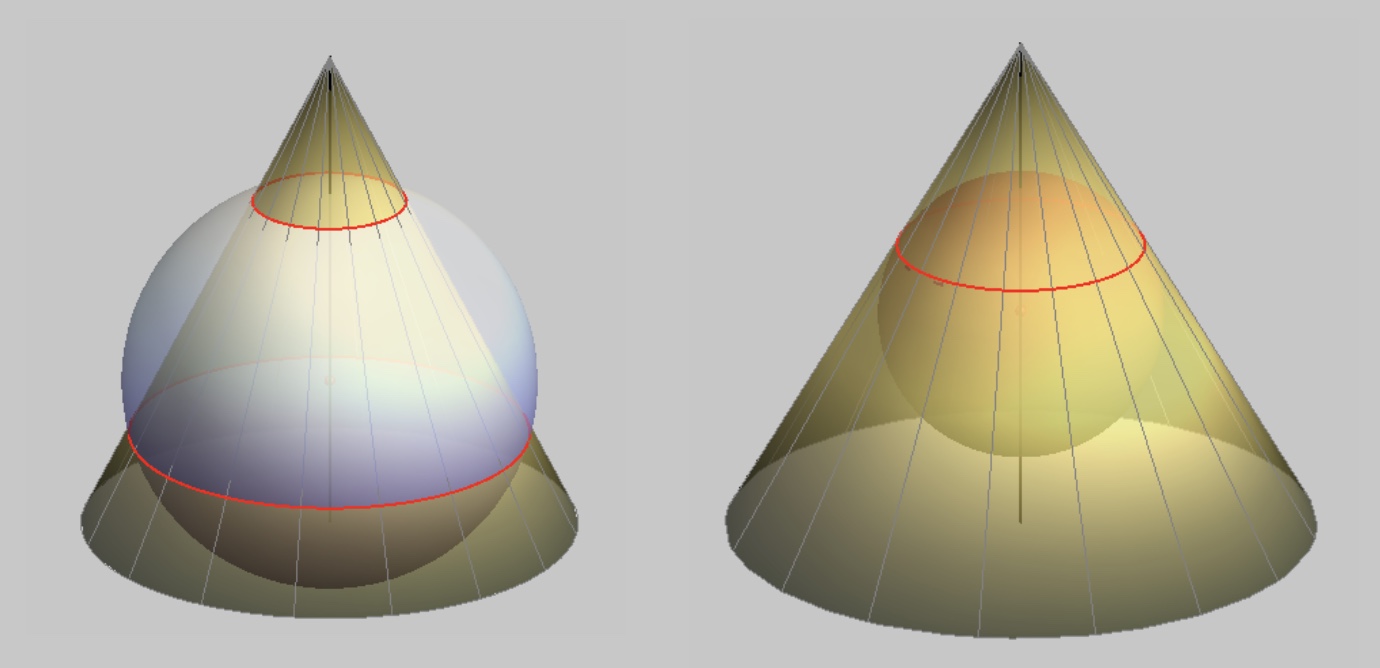} \ \ \ 
\includegraphics[scale = 0.1, trim = 0pt 0pt 700pt 0pt, clip]{ball-socps-more.JPG} \ \ \
\includegraphics[scale = 0.11]{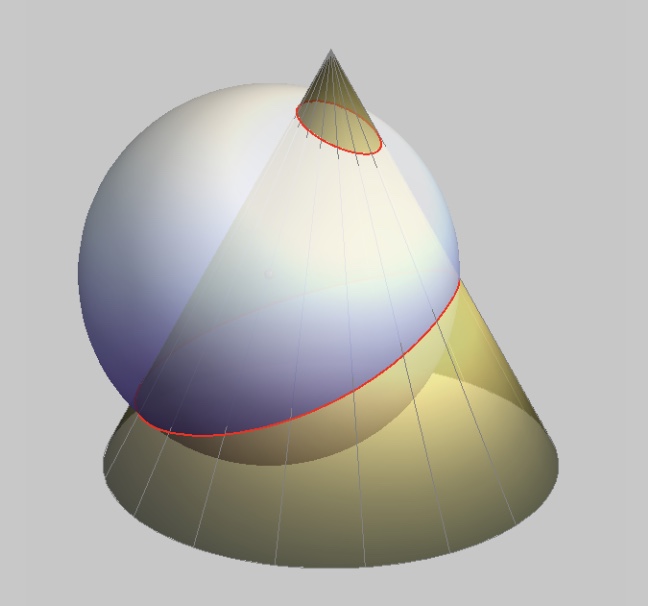}    \caption{Some further examples: \HC{s} may be trivial in the sense that one set is already contained in the other, or there may be parallel hyperplanes, or non-parallel hyperplanes. }
    \label{fig:cone-ball-example}
\end{figure}

\begin{figure}[h!]

\begin{center}
\begin{tikzpicture}

\definecolor{greenish}{RGB}{144,238,144}
\definecolor{blueish}{RGB}{173,216,230}
\coordinate (A) at (0,0);
\coordinate (B) at (-1,1);
\coordinate (D) at (-1,2);
\coordinate (C) at (-2,2);

\coordinate (G) at (1,1);
\coordinate (E) at (1,2);
\coordinate (F) at (2,2);
\coordinate (H) at (4.5,-0.5); %

\fill[blueish, opacity=0.3, draw=blue] (A) -- (B) -- (C) -- (D) -- (E) -- (G) -- cycle;
\fill[greenish, opacity=0.3, draw = green] (A) -- (B) -- (D)  -- (E) -- (F) -- cycle;

\draw[dashed, thick] (A) -- (B);
\draw[dashed, thick] (A) -- (G);
\draw[dashed, thick] (D) -- (E);

\node at (A) [below left] {A};
\node at (B) [below left] {B};
\node at (C) [above left] {C};
\node at (D) [above ] {D};

\node at (G) [below right] {G};
\node at (E) [above ] {E};
\node at (F) [above right] {F};

\end{tikzpicture}
\end{center}
\caption{Example where 3 hyperplanes are needed for a \HC, and all hyperplanes  intersect in an $n-1$ dimensional way.  
The convex sets are defined by their vertices ACEG and ABDF, respectively.
In this case, the union of the convex sets is preserved inside a cell of the hyperplane arrangement, but neither is contained in the other.}
\label{fig:n-1example}
\end{figure}
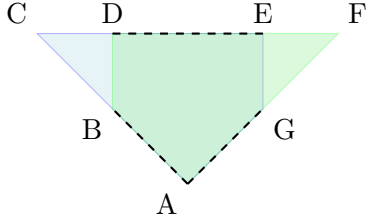

We show in Section~\ref{sec:idealBHC} that this leads to a simple structure in the hyperplane arrangement.
Unfortunately, not all convex sets have ideal boundary covers.  It may happen that  a cell of the hyperplane arrangement   might contain several convex sets that are weakly separable on the cell or, as we demonstrate, form a convex set in their union. See~\Cref{fig:types-bhc,fig:cone-ball-example} for examples.
Note that some of these complicating properties still persist even if all boundary intersections are $(n-1)$-dimensional.  See Figure~\ref{fig:n-1example}.

\subsection{Decomposition based on \HC} 
To deal with some of the complications from the prior examples, we prove the following key lemma in the generality where $K$ is a convex set.
However, we will only need this for $K$ being a polyhedron, as it represents a cell in our line of thought.
\begin{lemma}
\label{lem:boundary}
    Let $K$ be a convex set and let $C_1, C_2 \subseteq K$ be closed convex sets such that $\bd C_1 \cap \bd C_2 \subseteq \bd K$.  Suppose that $K$, $C_1$ and $C_2$ are full-dimensional. Then either 
    \begin{enumerate}
        \item $C_1$ and $C_2$ are weakly separable, or
        \item $C_1 \cup C_2$ is convex.
    \end{enumerate}
\end{lemma}
\begin{proof}
Suppose that $C_1$ and $C_2$ are not weakly separable.  Then the intersection $C_1 \cap C_2$ cannot be contained in a hyperplane.  Thus the intersection is full-dimensional and therefore 
there exists $\y \in C_1 \cap C_2$ and $\varepsilon > 0$ such that 
\begin{equation*} 
\label{eq:inner-y}
B(\y, \varepsilon) \subseteq C_1 \cap C_2.
\end{equation*}

We will prove that $C_1 \cup C_2$ is convex by contradiction.
So, assume $C_1 \cup C_2$ is not convex. 
Then, there exist $\x_1 \in C_1$ and $\x_2 \in C_2$ such that $\z_\lambda = \lambda \x_1 + (1-\lambda)\x_2 \notin C_1 \cup C_2$ for some $\lambda \in (0,1)$.
Without loss of generality, we can assume that $\x_i \in \partial C_i$ for $i = 1, 2$, and that 
 $\z_\lambda \notin C_1 \cup C_2$ for all $\lambda \in (0, 1)$.

Let $T = \conv(\{\y, \x_1, \x_2\})$.  Then by 
Theorem~\ref{thm:reverse_convex_polytope}, $W = \conv(T \setminus ( \intr(C_1)  \cup \intr(C_2)))$ is a polytope. We will look at the vertices of $W$.  Note that the half-open line segments $(\x_i,\y] \subset \intr(C_i)$ for $i=1,2$. These containments follow since $\y \in B(\y, \varepsilon) \subseteq \intr(C_i)$,  $\x_i \in C_i$, and $C_i$ is convex.  

Clearly $\x_1$ and $\x_2$ are vertices of $W$.  Since $C_1, C_2$ are closed and $\z_\lambda\notin C_1 \cup C_2$ there exists a $0 < \delta < 1$ such that  
$\z_{\lambda} + \delta (\y-\z_\lambda) \notin C_1 \cup C_2$.  This point is then also in $T$ and therefore in $W$ and not in $[\x_1, \x_2]$.  
Thus, we must have additional vertices of $W$.   Since the other edges of $T$ are contained in $C_1$ and $C_2$, any other vertex must lie  on the intersection of $\bd C_1 \cap \bd C_2$. 

Thus, there must be a point $\z^* \in \bd C_1 \cap \bd C_2 \subseteq \bd K$ with $\z^* \in \relint(T)$.  
However, $\relint(T) \subseteq \relint(K)$, and therefore $\z^* \notin \bd K$, which is a contradiction.
See Figure~\ref{fig:proof-lem-boundary}.

Therefore, $C_1 \cup C_2$ is convex if $C_1$ and $C_2$ are not weakly separable.
\end{proof}
\input{proof-picture}

In \Cref{sec:alternative-proof}, we provide a separate proof that only relies on the Intermediate Value Theorem and the continuity of the Minkowski-Gauge function.  However, this proof is more technical.

\begin{corollary}
\label{cor:boundary}
    Let $K$ be a convex set and let $C_1, \dots, C_m \subseteq K$ be convex sets such that $\bd C_i \cap \bd C_j \subseteq \bd K$ for every pair $(i,j)$ with $1 \leq i < j \leq m$.  Suppose that $K$, $C_1, \dots, C_m$ are full-dimensional. Suppose the pair $C_{[i]} := \cup_{s = 1}^i C_s$, $C_{i+1}$ is not weakly separable for all $i=1,\dots, m-1$. 
 Then $\bigcup_{i=1}^m C_i$ and $\bigcup_{i=1}^m \intr(C_i)$  are convex.
\end{corollary}
\begin{proof}
    We do this by induction.  If $m=2$, then the result follows directly from Lemma~\ref{lem:boundary}.  
    Suppose $m > 2$.  Then, using the inductive hypothesis, $C_{[m-1]}$
    is convex. 
    Lastly, the result follows by applying Lemma~\ref{lem:boundary} to $C_{[m-1]}$ and $C_m$. 
\end{proof}

\subsection{Decomposition based on Ideal Boundary Hyperplane Cover}

\label{sec:idealBHC}
We briefly delve into the Ideal \HC property.
This structure appears for Euclidean balls, as indicated above in~\Cref{subsec:BHC}.
\begin{proposition}
\label{prop:sphere-iBHC}
    Let $\{C_i\}_{i=1}^m$ be a collection of pairwise fully intersecting or disjunct Euclidean balls.
    Then $\{C_i\}_{i=1}^m$ has an Ideal \HC.  
\end{proposition}
\begin{proof}
    This follows from calculations in Proposition~\ref{lem:spheres}.
\end{proof}

\begin{lemma}[Ideal Hyperplane Separation]
\label{lem:ideal-separation}
Let $C_{1}, C_{2}$ be full-dimensional closed convex subsets of $\mathbb{R}^{n}$ with an ideal boundary cover $H$ (see Definition~\ref{defn:ideal-boundary-cover}).

Then, either 
$$
C_1 \cap H^< \subseteq C_2 \cap H^< \qquad \text{or} \qquad  C_2 \cap H^< \subseteq C_1 \cap H^<,
$$
where $H^< := \intr(H^{\leq})$.
\end{lemma}

\begin{proof}
First, if for either $i=1,2$ we have $C_i \cap H^< = \emptyset$, then the result is trivial.

So assume that there exist $\x^i \in C_i \cap H^<$ for $i=1,2$.  

Note that $H^{\leq}$ is convex and $C_i \cap H^{\leq} \subseteq H^{\leq}$, $i = 1, 2$, are closed convex. 
Even further, $\partial (C_1 \cap H^{\leq}) \cap \partial (C_2 \cap H^{\leq}) = (\partial C_1 \cap \partial C_2) \cap \partial H^{\leq} \subseteq H$, as $H = \partial H^{\leq}$ and $H$ represents the ideal cover.
Hence, we can apply~\Cref{lem:boundary} to $C_1 \cap H^{\leq}$, $C_2 \cap H^{\leq}$, and $H^{\leq}$.

It follows that either $C_1 \cap H^{\leq}$ and $C_2 \cap H^{\leq}$ are weakly separable or 
$$
(C_1 \cap H^{\leq}) \cup (C_2 \cap H^{\leq}) = (C_1 \cup C_2) \cap H^{\leq}
$$ 
is convex.

If we are in the weakly separable case, there exists a hyperplane $H'$ which contains $\partial(C_1 \cap H^{\leq}) \cap \partial(C_2 \cap H^{\leq})$ and separates $\x^1$ and $\x^2$.
As $\dim(\aff(\partial C_1 \cap \partial C_2)) = n - 1$ due to the ideal boundary cover and
$$ 
\partial C_1 \cap \partial C_2 \subseteq (C_1 \cap H^{\leq}) \cap (C_2 \cap H^{\leq}) \subseteq H',
$$
we can take the affine hull to get
$$
H = \aff(\partial C_1 \cap \partial C_2) \subseteq H'.
$$
Thus, $H = H'$ since $H$, $H'$ are $(n-1)$-dimensional affine spaces. 
However, this creates a contradiction since $H'$ is assumed to separate $\x^1$ and $\x^2$, which is impossible since $\x^1, \x^2 \in H^<$, i.e., are not separated by $H$,  and $H = H'$.

Therefore, $(C_1 \cup C_2) \cap H^{\leq}$ is convex and so is $\intr((C_1 \cup C_2) \cap H^<)$ which is equivalent to $\intr(C_1 \cup C_2) \cap H^<$, see~\Cref{obs:intr-intersection}.
In other terms, $C_1 \cup C_2$ is convex in the open subspace $H^<$.
From above, we know that $\partial C_1 \cap \partial C_2 \subseteq H$, thus the boundaries of $C_1$ and $C_2$ do not intersect in $H^<$.
Together, this gives that in $H^<$, we have $C_1 \subseteq C_2$ or vice versa.
This shows the claim.

\end{proof}

We now show why this property is useful.

\begin{theorem}
\label{thm:ideal-decomposition}
    Let $\{C_i\}_{i=1}^m$ be closed, full-dimensional convex sets with an Ideal \HC $\mathcal H$.  
    Then, on every maximal cell of the corresponding hyperplane arrangement, for any $C_i, C_j$ that intersect the cell on its interior, we have that either $C_i$ and $C_j$ are strongly separable on the cell, or one set contains the other when restricted to the cell.
\end{theorem}
\begin{proof}
    This follows from Lemma~\ref{lem:ideal-separation}.
\end{proof}

\subsection{Decomposition based on General \HC}
We now state our first main structural result for decomposing unions of convex sets admitting a boundary hyperplane cover. The construction assumes access to exact separation oracles, which may be unavailable in practice. In~\Cref{sec:main-integer} we present a tractable variant that, invoking \Cref{lem:separation-integer-hulls}, separates the integer points of the sets rather than the sets themselves, and thus guarantees decomposition only on the integer lattice. This formulation enables the subsequent analysis of algorithmic complexity. The proofs follow parallel outlines, but the critical ingredients—and hence the implications—differ.

\begin{figure}
    \centering

\begin{center}
\includegraphics[scale = 0.18]{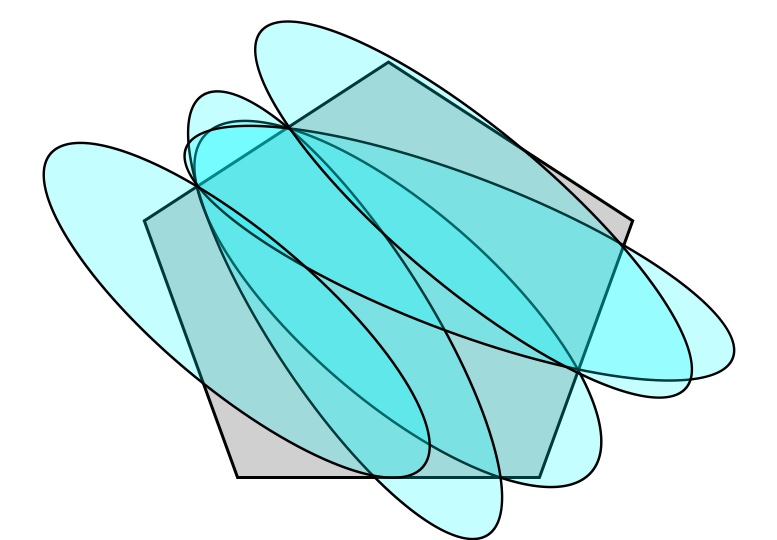}
\includegraphics[scale=0.18,trim=0cm 0.02cm 0cm 0.015cm,clip]{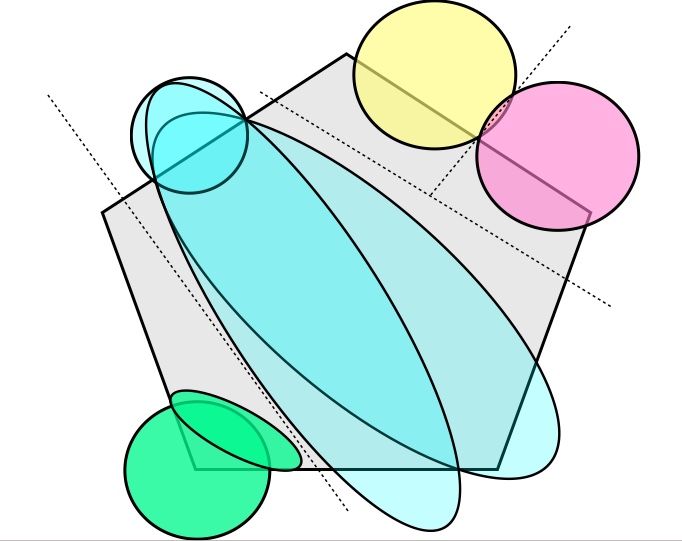}
\includegraphics[scale=0.18,trim=0cm 0.05cm 0cm 0.07cm,clip]{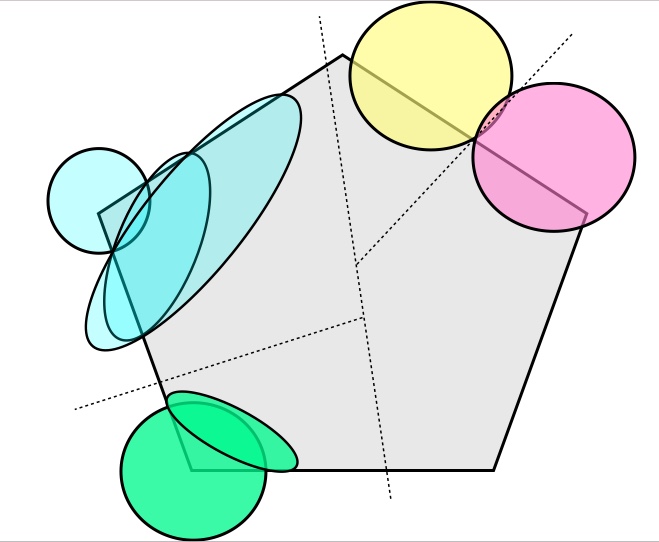}

\end{center}
    \caption{Various cases of how the convex sets can interact on a maximal cell.  Connected components in the intersection graph are color-coordinated.}
    \label{fig:cell-decomposition}
\end{figure}

\begin{theorem}[Decomposition given \HC (exact separation, i.e, continuous version)]
\label{thm:main}
Let $\mathcal C$ be a class of bounded convex sets that is closed under polyhedral intersections and satisfy the inradius intersection property (\Cref{assumption:intersection-radius}) and suppose we are provided the oracle $\WEAK{\mathcal C}{\R^n}$ (Oracle~\ref{oracle:separation}).  Let  $C_1, \dots, C_m\in \mathcal C$ be full-dimensional closed and bounded convex sets with a \HC~$\mathcal H$.
Then there exist convex sets $K_s$ with pairwise disjoint interiors such that 
\begin{equation*}
    \bigcup_{i=1}^m \intr(C_i) = \bigcup_{s=1}^t K_s.
\end{equation*}
Specifically, there exist polyhedra $P_s$ with pairwise disjoint interiors such that $K_s = P_s \cap \bigcup_{i \in I_s} \intr(C_i)$, where $I_s \subseteq [m]$, and
\begin{equation*}
\R^n \setminus \bigcup_{i=1}^m \intr(C_i)
 = \bigcup_{s=1}^t (P_s \setminus K_s).
\end{equation*}
Here,  $t = O((m^2d)^n)$ where $d = |\mathcal H|$. 
Note that the set differences $P_s \setminus K_s$, having again pairwise disjoint interiors, form a decomposition of $\R^n \setminus \bigcup_{i = 1}^t\intr(C_i)$.
Each component of this decomposition is then defined by the triple $(P_s, K_s, I_s)$.

We can compute such a decomposition in oracle-polynomial time with the dimension $n$ fixed.
\end{theorem}
Observe that the sets $K_s$ might be half open.  
\begin{proof}
Enumerate the maximal cells $\mathcal P_{\max}$ of the hyperplane arrangement $\mathcal H$, of which there are at most $|\mathcal P_{\max}| = O(d^n)$ by Lemmas~\ref{lem:hyperplane-arrangement-bound},~\ref{lem:hyperplane-arrangement}. 
For each $P \in \mathcal P_{\max}$, let $I_P \subseteq [m]$ be the set of indices $i$, for which convex set $\intr(C_i)$ intersects $P$.  
As $\intr(C_i)$ and $P$ are full-dimensional, and $\intr(C_i)$ is an open set, their intersection is either empty or full-dimensional.
In particular, we use~\Cref{cor:intersection-full} combined with~\Cref{lem:interior-intersection-full-dim} to detect this.

\noindent \underline{Intermediate step.}
We claim that 
\begin{equation}
\label{eq:polytope-decomposition}
\R^n \setminus \bigcup_{i=1}^m \intr(C_i)
 = \bigcup_{P \in \mathcal P_{\max}} (P \setminus \bigcup_{i \in I_P} \intr(C_i)).
\end{equation}

To prove this, let \( \x \in \R^n \). Since \( \{P\}_{P \in \mathcal P_{\max}} \) tiles \( \R^n \), there exists at least one \( P \in \mathcal P_{\max} \) such that \( \x \in P \). If \( \x \) lies on the boundary between multiple polytopes, then there exist maximal cells \( P_1, \dots, P_q \in \mathcal P_{\max} \) such that \( \x \in \partial P_j \) for each \( j = 1, \dots, q \).
We now proceed with the analysis of \( \x \) with respect to the sets \( C_i \) and the polytopes \( P_1, \dots, P_q \).

\textbf{Case 1:} \( \x \notin \bigcup_{i=1}^m \intr(C_i). \)

In this case, the point \( \x \) is outside the interiors of all the sets \( C_i \). Therefore, \( \x \) will appear on both the left-hand side (LHS) and in all the relevant terms on the right-hand side (RHS) of \eqref{eq:polytope-decomposition} for the maximal cells \( P_1, \dots, P_q \) that contain \( \x \).
This confirms that the decomposition holds in this scenario.

\textbf{Case 2:} \( \x \in \intr(C_i) \) for some \( i \in \{1, \dots, m\}. \)

If \( \x \) lies in the interior of some \( C_i \), then \( \x \) is excluded from both sides of \eqref{eq:polytope-decomposition}. Specifically, since \( C_i \in I_{P_j} \) for all \( j \in \{1, \dots, q\} \), \( \x \) does not appear on the LHS or in any of the terms on the RHS. Thus, the equation holds in this case as well.

Since both cases hold, we conclude that \eqref{eq:polytope-decomposition} is valid.\\

\noindent \underline{Restrict to a single cell.}
We now restrict our attention to some $P\in \mathcal P_{\max}$, and without loss of generality, $I_P = \{1, \dots, m'\}$ for some $m' \leq m$. 
Note that for any $C_i$ and $C_j$ with $i,\, j \in I_P$, $i \neq j$, it holds true that the possibly empty set $\partial C_i \cap \partial C_j \cap P$ is a subset of $\partial P$ according to the construction of $P$ via the \HC property.
In particular, we get that 
\begin{equation}
	\label{eq:boundary-intersect-on-cell}
	\partial(C_i \cap P) \cap \partial(C_j \cap P) \subseteq \partial P
\end{equation}

\textbf{Intersection detection.}  
We build a graph $G = ([m'], E)$ with nodes $[m']$ and edges $(i,j) \in E$ whenever $\intr(C_i) \cap \intr(C_j) \cap P \neq \emptyset$, in other terms, the full-dimensional sets $C_i$ and $C_j$ intersect in~$P$ and are not weakly separable. 
Algorithmically, this requires access to the oracle $\WEAK{\mathcal C}{\R^n}$ (Oracle~\ref{oracle:separation}).  In particular, we use \Cref{rem:oracle-C-intersect-P} or \Cref{thm:full-dim} combined with \Cref{obs:intr-intersection} to detect this.
From above we know that $\intr(C_i \cap P) \neq \emptyset$ and $\intr(C_j \cap P) \neq \emptyset$.
Further, \Cref{assumption:intersection-radius} transfers to $C_i \cap P$ and $C_j \cap P$.
Together with the boundedness of the convex sets this allows to apply~\Cref{cor:intersection-full} for $C_i \cap P$ and $C_j \cap P$.
This returns whether $\intr(C_i \cap C_j \cap P)$ is empty and with \Cref{obs:intr-intersection} and \Cref{lem:interior-intersection-full-dim} (where $``C_2" = P$) we can simultaneously conclude whether $\intr(C_i) \cap \intr(C_j) \cap P$ is empty.

a.) If $G$ is connected, then 
we can produce a \emph{connected ordering} $[i_1, \dots, i_{m'}]$ such that each of the  induced subgraphs $G_j$ with nodes $[i_1, \dots, i_j]$ are connected for $j=1, \dots, m'$.  For instance this can be achieved using a breadth-first-search in time $O(m' + |E|) = O(m'^2)$. Given the new ordering, we have 
$$
\left(\bigcup_{i = 1}^j (C_i \cap P)\right) \cap (C_{j+1} \cap P) = P \cap \bigcup_{i = 1}^j C_i \cap C_{j+1} \neq \emptyset
$$ 
for $j \in [m' - 1]$. In addition, we have that $\bigcup_{i = 1}^j C_i \cap P$ and $C_{j+1} \cap P$ are not weakly separable for $j \in [m' - 1]$.

Together with~\eqref{eq:boundary-intersect-on-cell}, this allows to apply Corollary~\ref{cor:boundary} for $P$ and its convex subsets $C_i \cap P$, $i \in [m']$.
It follows that $\bigcup_{i = 1}^{m'} \intr(P \cap C_i) = \intr(P) \cap \bigcup_{i=1}^{m'} \intr(C_i)$ is convex. 
In particular, using
~\Cref{lem:interior-intersection-full-dim}, 
$K := P \cap \bigcup_{i = 1}^{m'} \intr(C_i)$ is convex and we can return $(K,\,P,\, [m'])$ as one of our components in the decomposition.

b.) Alternatively, if $G$ is disconnected, then on each connected component, we can make a similar argument that each connected component forms a convex set.  Let $K_1, \dots, K_{m''}$ be these convex sets with $m'' \leq m' \leq m$ where indices $I_v \subseteq [m']$ are such that $P \cap C_i \subseteq K_v$ for $i \in I_v, v \in [m'']$.
In particular,
$$
K_v = P \cap \bigcup_{i \in I_v} \intr(C_i).
$$
For $P$ and the tuples $(K_v, I_v)$, $v \in [m'']$, we apply the next step to compute corresponding polyhedra $P_v$ that are disjoint and intersect on the respective $K_v$ only.

\textbf{Separation.} 
For this, using $\WEAK{\mathcal C}{\R^n}$  and the outer representation of $P$, we compute a weak separating hyperplane $H_{\a^{vw},b_{vw}}$ of the closures of $K_v$ and $K_w$ for any pair $v,\,w \in [m'']$, $v < w$.  \Cref{rem:oracle-union} provides the separation oracles requested in the definition of $\WEAK{\mathcal C}{\R^n}$.
Note that these are ${ m'' \choose 2}$ many pairs to consider.

Without loss of generality, we orientate the corresponding halfspaces such that $K_v \subseteq H_{\a^{vw},b_{vw}}^{\leq}$ if $v < w$ and $K_v \subseteq H_{\a^{wv},b_{wv}}^{\geq}$ if $w < v$.
Now, we set 
\begin{equation}
\label{eq:Pv-def}
P_v := P \cap \bigcap_{v < w} H_{\a^{vw},b_{vw}}^{\leq} \cap \bigcap_{w < v} H_{\a^{wv},b_{wv}}^{\geq},
\end{equation}
for all $v \in [m'']$, where $w$ iterates over $[m'']\setminus\{v\}$.
We then return $(K_v,\, P_v,\, I_v)$ as components in the decomposition. 

As there remain polyhedra in the scheme of the right-side of~\labelcref{eq:Pv-def} when following all the other $v-w$-combinations, we detect all such non-empty ones by means of~\Cref{lem:hyperplane-arrangement}, say $P'_1,\dots,P'_{s'}$.
We can return them with $K'_1 = \dots = K'_{s'} = \emptyset$ and index sets $I'_1 = \dots = I'_{s'} = \emptyset$.

Since the amount of hyperplanes required for the definition all polyhedra in this paragraph is $m''(m''-1)/2$, as noted above, by~\Cref{lem:hyperplane-arrangement-bound}, we have that $m'' + s' = O(m^{2n})$.
Since we fixed a cell $P$ and there are at most $O(d^n)$ such, we receive $O((m^2d)^n)$ many pieces. 
Then we have that 
\begin{equation}
\label{eq:P-decomposition}
    P \setminus \bigcup_{i \in I_P} \intr(C_i) = \left(\bigcup_{s = 1, \dots ,m''}   P_s \setminus K_s \right) \cup \left(\bigcup_{s = 1, \dots ,s'}   P'_s \setminus K'_s\right).
\end{equation}

Combining \eqref{eq:polytope-decomposition} and \eqref{eq:P-decomposition} completes the result.
\end{proof}

\begin{figure}
	\centering
	\begin{tikzpicture}
		
		\node (A) at (0, 0) {\input{heart.tex}};
		\node (C) at (5.2, -1.2) {\input{heart1.tex}};
		\node (B) at (5.2, 1.2) {\input{heart2.tex}};
		
		\draw[->, thick] (A.east) -- (B.west);
		\draw[->, thick] (A.east) -- (C.west);
		
		\node (D) at (10, 1.2) {\input{heart2a.tex}};
		\node (E) at (10, -1) {\input{heart2b.tex}};
		\draw[->, thick] (B.east) -- (D.west);
		\draw[->, thick] (B.east) -- (E.west);
		
	\end{tikzpicture}
	\caption{An example of the procedure in Theorem~\ref{thm:main}.  This example initially has two non-empty cells of the hyperplane arrangement.  The triangle cell demonstrates a case where the union of convex sets on that cell is a larger convex set.  The other non-empty cell shows that the convex sets are now separable on this cell, so we further decompose the region into two sub-cells.  The leaves of this tree represent all the pieces that we need.}
	\label{fig:heart-decomposition}
\end{figure}

See~\Cref{fig:heart-decomposition} for a simple example of this procedure
and  Figure~\ref{fig:cell-decomposition} for examples of how a maximal cell can appear.

Note that the decomposition provided in Theorem~\ref{thm:main} is not unique.  In fact, it is likely more efficient to apply separation techniques ahead of the decomposition and handle contiguous unions of convex sets more individually.  However, we find that the approach in the proof is cleaner to present.

Alternatively,  the above theorem is much simpler to prove in the context of an Ideal \HC, following~\Cref{thm:ideal-decomposition}. Having such an Ideal \HC indeed would be valuable to consider to make such a decomposition more practical and easier to achieve.

In the absence of numerically stable exact separation of some convex sets, we prove an alternative decomposition in \Cref{thm:main-integer} that allows for precision to be taken, since we care most about feasibility over integer points.  For this result, we build up results that are based on discrete points and also define a new separation oracle that applies to the interior of a convex set.

\section{Decomposition of Integers via \HC }
\label{sec:main-integer}

We now seek to avoid numerical issues in the prior results by focusing on integer points.  For example, this will allow us to avoid the assumption of a known inradius.

We will need some formality to work with interiors of convex sets.
Perhaps non-standard, we define a separation oracle over the interior of a convex set.

\begin{oracle}[Interior Separation Oracle]
\label{oracle:interior-separation}
An \emph{interior separation oracle}  \( \SEPint{C}{D} \) for a convex set $C \subseteq \R^n$ determines whether a given point from a domain $D$ is in the interior of $C$ or provides a hyperplane separating the point from the interior of $C$.

\end{oracle}

A key distinction of this second oracle is that it returns a valid inequality even when \( \bar{\x} \in \partial C \). This assumption is natural in many settings and aligns with the principle of returning inequalities that are tight for \( C \). For example, if \( C \) is a polyhedron, the oracle returns tight valid inequalities. If instead \( C = \{\x \in \mathbb{R}^n : f(\x) \leq 0\} \) for some convex function \( f \), then the oracle provides subgradients of \( f \), even at boundary points where \( f(\bar{\x}) = 0 \). See, for instance, \cite{OERTEL2014424} for an integer programming algorithm that makes use of such subgradients.

\subsection{Integer Feasibility  on the Interior of a Convex Set}

In some cases, we may want to check integer feasibility on the interior of a convex set.  
We show that we can adapt a separation oracle to apply to a new convex set that has no integer points on its boundary, and hence has the same integer points on its interior.
We begin with a simple lemma.

\begin{lemma}[Shifted Inequality Cuts Off Boundary Integer Points]
\label{lem:shifted-separation-integral}
Let \( C \subseteq \mathbb{R}^n \) be a convex set 
and suppose \( \a^\top \x \leq b \) is a valid inequality for \( C \), with \( \a \in \mathbb{Z}^n \), \( b \in \mathbb{Z} \), and \( \a \neq \ve 0 \). Then the \emph{shifted inequality}
    \(
    \a^\top \x \leq b - 1
    \)
    is valid for all integer points in \( \operatorname{int}(C) \cap \mathbb{Z}^n \).
\end{lemma}

\begin{proof}
Since \( \a^\top \x \leq b \) is valid for all \( \x \in C \), and \( \operatorname{int}(C) \subseteq C \), we have
\(
\a^\top \x < b \quad \text{for all } \x \in \operatorname{int}(C).
\)
Now let \( \x \in \operatorname{int}(C) \cap \mathbb{Z}^n \). Because \( \a^\top \x \in \mathbb{Z} \) and \( \a^\top \x < b \), it follows that
\(
\a^\top \x \leq b - 1.
\)
Therefore, the shifted inequality \( \a^\top \x \leq b - 1 \) is valid for all \( \x \in \operatorname{int}(C) \cap \mathbb{Z}^n \).
\end{proof}

\begin{lemma}[Separation for Interior Integer Hull via Interior Oracle]
\label{lem:integer-hull-from-interior-sep}
Let \( C \subseteq \R^n \) be a convex set with \( C \subseteq B_\infty(\ve 0, R) \), and let \(\SEPint{C}{\Z^n}\) be the interior separation oracle for \(C\). Define
\[
C' := \Biggl( \bigcap_{\substack{\x \in \Z^n \cap B_\infty(\ve 0,R) \\ \x \notin \intr(C)}} 
\left\{ \z \in \R^n : \a_\x^\top \z \leq b_\x - 1 \right\} \Biggr) \cap B_\infty(\ve 0,R),
\]
where for each such \(\x\), the inequality \((\a_\x,b_\x)\) is returned by \(\SEPint{C}{\Z^n}(\x)\).  

Then the polytope \(C'\) satisfies
\[
C' \cap \Z^n \;=\; \intr(C) \cap \Z^n,
\]
and moreover admits a separation oracle \(\SEP{C'}{\Z^n}\) whose complexity is bounded as
\[
\Phi_{\SEP{C'}{\Z^n}}(\tau) = \operatorname{poly}\!\bigl(\Phi_{\SEPint{C}{\Z^n}}(\tau), \langle R \rangle \bigr).
\]
\end{lemma}

\begin{proof}
First we show correctness of the integer hull equality.  

($\subseteq$) Suppose \(\x \in \intr(C) \cap \Z^n\). By Lemma~\ref{lem:shifted-separation-integral}, for $\y \notin \intr(C)$, every shifted inequality \(\a_\y^\top z \leq b_\y - 1\) is valid for \(\x\). 
Thus \(\x \in C'\).

($\supseteq$) Conversely, let \(\x \notin \intr(C)\cap \Z^n\). 
Then, either $\x \notin B_\infty(\ve 0, R)$, or calling $\SEPint{C}{\Z^n}$ on $\x$ returned an inequality $(\a_\x, b_\x)$ in the definition of $C'$. 
In the former case, it directly follows that $\x \notin C'$. 
In the latter case, it holds true that $\a^\top \x = b_\x$ and, thus, $\x \notin C'$ according to the definition of $C'$.

To implement the oracle \(\SEP{C'}{\Z^n}\) on input \(\bar{\x} \in \Z^n \cap B_\infty(\ve 0,R)\) first query $\SEPint{C}{\Z^n}$ on $\bar{\x}$.
\begin{itemize}
    \item If it confirms membership in $\intr(C)$, confirm it in \(C'\).  
    \item Otherwise, consider the returned \((\a_{\bar{\x}}, b_{\bar{\x}})\) and output the shifted inequality \(\a_{\bar{\x}}^\top \z \leq b_{\bar{\x}} - 1\) as a separating hyperplane.  
\end{itemize}
In addition, membership in \(B_\infty(\ve 0,R)\) is enforced by its \(2n\) inequalities, which can be checked explicitly.  

Thus \(\SEP{C'}{\Z^n}\) is implementable with polynomial overhead in the complexity of \(\SEPint{C}{\Z^n}\) and \(\langle R \rangle\), as claimed.
\end{proof}

Using this lemma we can optimize over the  the interior integer points of a convex set.  To the best of our knowledge, such a result has never been stated.
\begin{theorem}[Interior Convex Integer Feasibility and Optimization]
\label{thm:interior-optimization}
Let \( C \subseteq \R^n \) be a convex set with \( C \subseteq B_\infty(\ve 0,R) \), 
and suppose we have access to the interior separation oracle \(\SEPint{C}{\Z^n}\). 
Then there exist algorithms which

\begin{itemize}
  \item either find a point in \( \intr(C) \cap \Z^n \) or correctly decide that 
  \( \intr(C) \cap \Z^n = \emptyset\), and
  \item given \( \mathbf{c} \in \Z^n \), solve
  \(
  \max \{ \mathbf{c}^\top \x : \x \in \intr(C) \cap \Z^n \}.
  \)
\end{itemize}

The algorithm runs in time
\(
2^{O(n \log n)} \, \operatorname{poly}\!\bigl(\langle R \rangle, \Phi_{\SEPint{C}{\Z^n}}(n \langle R\rangle )\bigr).
\)
\end{theorem}

\begin{proof}
By Lemma~\ref{lem:integer-hull-from-interior-sep}, from 
\(\SEPint{C}{\Z^n}\) we can construct a convex set \(C' \subseteq \R^n\) with
\[
C' \cap \Z^n = \intr(C) \cap \Z^n,
\]
together with a separation oracle \(\SEP{C'}{\Z^n}\) whose complexity is polynomial in 
\(\Phi_{\SEPint{C}{\Z^n}}(\cdot)\) and \(\langle R \rangle\). 
Applying Theorem~\ref{thm:convex-integer} to \(C'\) with the oracle 
\(\SEP{C'}{\Z^n}\) gives the feasibility and optimization algorithms with complexity
\[
2^{O(n \log n)} \, \operatorname{poly}\!\bigl(\langle R \rangle, \Phi_{\SEP{C'}{\Z^n}} (n \langle R \rangle) \bigr).
\]
Since \(\Phi_{\SEP{C'}{\Z^n}}(\cdot)  = \operatorname{poly}(\Phi_{\SEPint{C}{\Z^n}} (\cdot) , \langle R \rangle)\),
the claimed complexity bound follows.
\end{proof}

\subsection{Weak Separation Over the Integers}

We now consider the weak separation oracle (\Cref{oracle:separation}) in the context of separating integer points.  
This will be easier to be numerically precise.

\begin{lemma}[Weak Separation of Integer Points]
\label{lem:separation-integer-hulls}
Let $R \in \Z_+$ and let $\mathcal C$ be a collection of (possibly open) convex sets contained in $B(\ve 0, R)$ and equipped with oracles to optimize of a linear function over the integer points in the convex sets.
Suppose that $\mathcal C$ is closed under finite intersections with halfspaces.

Then there is an algorithm for $\WEAK{\mathcal C}{\Z^n}$ that runs in oracle-polynomial time in fixed-dimension.

In particular, for two open convex sets $C_1,\, C_2 \in \mathcal C$ the presented with interior separation oracles $\SEPint{C_i}{\Z^n}$ there is an algorithm to compute $\WEAK{\mathcal C}{\Z^n}(C_1, C_2)$ that runs  in 
\[
n2^{O(n \log{n})} \operatorname{poly}(\langle R \rangle, \Phi_{\SEPint{C_1}{\Z^n}}(n\langle R \rangle), \Phi_{\SEPint{C_2}{\Z^n}}(n\langle R \rangle)).
\]
\end{lemma}

\begin{proof}

 First we can check if $C_1$ and $C_2$ contain integer points using Theorem~\ref{thm:convex-integer} with the separation oracles of $C_1$ and $C_2$.  If either does not, then we return the inequality $0 \leq 1$.
 
Assuming that they are both containing integer points, we now look for a hyperplane $H_{\a,b}$ that weakly separates the points $C_1 \cap \Z^n$ and $C_2 \cap \Z^n$. 
Thus, we seek to find a pair $(\a,b)$ with $\a \not\equiv 0$ such that 
\begin{equation}
   \label{eq:cuts} 
\begin{aligned}
    & \a^\top \bar{\x} \leq b, \quad \text{for all }  \bar \x \in C_1 \cap \mathbb{Z}^n, \\
& \a^\top \bar{\y} \geq b, \quad \text{for all }\bar \y \in C_2 \cap \mathbb{Z}^n.
\end{aligned}
\end{equation}

If  these points are indeed weakly separable, such a non-zero $\a$ must exist.  
As the separation property of $(\a, b)$ is invariant under positive scaling, there exists a feasible integral solution to at least one of the $2n$ subproblems where, for each $i$, either  
    $a_i \geq 1$ or $a_i \leq -1.$
In particular, let 
$$
P^{\pm}_i = \{ (\a,b) \in \R^{n+1} :  \eqref{eq:cuts} \text{ holds},\, \pm a_i \geq 1\},
$$
and for each of the $2n$  polyhedra $P^\pm_i$, decide the non-emptiness of each polyhedron; we show how this decision can be carried out efficiently.  
\paragraph{Checking non-emptiness of $P^{\pm}_i$.}
If $P^{\pm}_i \neq \emptyset$, by \cite[Corollary 17.1b]{schrijver_theory_1986} (which establishes polynomial correspondence between encoding size of vertices and facets of a rational polyhedron), there exists an integral solution of encoding size polynomially bounded in the encoding size of $\bar \x$ and $\bar \y$ from \eqref{eq:cuts}.  
In particular, there is a solution $(\a,b) \in \Z^n \times \Z$ with $\|(\a,b)\|_\infty \leq U$ where $\langle U \rangle$ depends polynomially on $\langle R \rangle$.

Now, we will show that a separation oracle can be implemented for a bounded version of $P^\pm_i$ and thus employ the optimization oracles 
Theorem~\ref{thm:convex-integer} to find a feasible integer solution.  
This prevents the need from explicitly writing down all of the inequalities that describe $P^\pm_i$.
We handle only the separation over $P^+_i$ as $P^-_i$ is similar.  
See~\Cref{alg:separation} of how to do this.

\begin{algorithm}
\caption{Separation Oracle $\SEP{P_i^+}{\Z^n}$ for $P^{+}_i$}
\label{alg:separation}
\begin{algorithmic}[1]
\REQUIRE A pair $(\a, b)$, sets $C_1$, $C_2$ presented with separation oracles
\ENSURE A valid inequality violated by $(\a, b)$ or a certificate of feasibility

\IF{$a_i \not\geq 1$}
    \RETURN the violated inequality  $a_i \geq 1$
\ENDIF

\STATE Solve $\bar{\x} \leftarrow \arg\max \{ \a^\top \x : \x \in C_1 \cap \mathbb{Z}^n \}$
\IF{$\a^\top \bar{\x} > b$}
    \RETURN inequality $\a^\top \bar{\x} \leq b$
\ENDIF

\STATE Solve $\bar{\y} \leftarrow \arg\min \{ \a^\top \y : \y \in C_2 \cap \mathbb{Z}^n \}$
\IF{$\a^\top \bar{\y} < b$}
    \RETURN inequality $\a^\top \bar{\y} \geq b$
\ENDIF

\RETURN Certify that $(\a, b)$ is feasible for $P_i^+$
\end{algorithmic}
\end{algorithm}

The argmin and argmax are  computed via their optimization oracles.  If we assume the sets are presented with interior separation oracles, then  \Cref{thm:interior-optimization} with the separation oracles of $C_1$ and $C_2$. 
In particular, $\SEP{P_i^+}{\Z^n}$ is encoded in time
\begin{equation*} 
	2^{O(n\log(n))}\operatorname{poly}(\langle R \rangle, \Phi_{\SEPint{C_1}{\Z^n}}(n\langle R \rangle ),\Phi_{\SEPint{C_2}{\Z^n}}(n\langle R \rangle ) ).
\end{equation*}
 Since $C_1$ and $C_2$ are bounded, the encoding size of the returned inequalities are bounded.

\paragraph{Conclusion.}
Finding a feasible integer point \((\a, b)\) using Theorem~\ref{thm:convex-integer} together with the above separation oracle---which itself internally uses Theorem~\ref{thm:convex-integer}---has total complexity
\begin{align*}
& n2^{O(n \log{n})} \operatorname{poly}(\langle R \rangle, \Phi_{\SEP{P_i^+}{\Z^n}}(n\langle R \rangle)) \\
= \ & n2^{O(n \log{n})} \operatorname{poly}(\langle R \rangle, \Phi_{\SEPint{C_1}{\Z^n}}(n\langle R \rangle), \Phi_{\SEPint{C_2}{\Z^n}}(n\langle R \rangle)).
\end{align*}
The additional factor $n$ stems from the $2n$ polyehdra which are to check.

If all $P^\pm_i$ are shown to not have an integer point, then we conclude that $C^I_1$ and $C^I_2$ cannot be weakly separated.

\end{proof}
Note that of course we did not need to seek an integral solution for $(\a,b)$.  
We of course could have just found a rational solution.  
This would result in a better time complexity using some variant of the ellipsoid algorithm.  However, this improvement is unimportant to the overall goal of this paper.

\subsection{Decomposition over the integers}
We're now prepared to present our main theorem about decomposition over the integers.

\begin{theorem}[Decomposition given \HC (integer version)]
\label{thm:main-integer}
    Suppose $C_1, \dots, C_m\subseteq \R^n$ are full-dimensional closed and bounded convex sets contained in $B(\ve 0, R)$  and have a \HC~$\mathcal H$.
Then there exist convex sets $K_s$ with pairwise disjoint interiors such that 
\begin{equation*}
    \Z^n \cap \bigcup_{i=1}^m \intr(C_i)   = \Z^n \cap \bigcup_{s=1}^t K_s.
\end{equation*}
Specifically, there exist polyhedra $P_s$ with pairwise disjoint interiors such that $K_s = P_s \cap \bigcup_{i \in I_s} \intr(C_i)$, where $I_s \subseteq [m]$, and
\begin{equation*}
    \Z^n \setminus \bigcup_{i=1}^m \intr(C_i) = \bigcup_{s=1}^t (P_s \cap \Z^n) \setminus K_s.
\end{equation*}
Here,  $t = O((m^2d)^n)$ where $d = |\mathcal H|$.
Note that the set differences $(P_s \cap \Z^n) \setminus K_s$ form a decomposition of $\Z^n \setminus \bigcup_{i = 1}^m \intr(C_i)$. 
Each component of the decomposition is then defined by the triple $(P_s, K_s, I_s)$.

\textbf{Complexity.} Suppose that $C_i$ are presented  by interior separation oracles $\SEPint{C_i}{\Z^n}$.
Then in  polynomial time with the dimension $n$ fixed, we can compute this decomposition in time polynomial in $m$, $d^n$, the encoding size $\phi$ of the hyperplanes in $\mathcal H$, and $\Phi_{\SEPint{C_i}{ \Z^n}}$ for $i=1, \dots, m$.
\end{theorem}

\begin{proof}
	Enumerate the maximal cells $\mathcal P_{\max}$ of the hyperplane arrangement $\mathcal H$, of which there are at most $|\mathcal P_{\max}| = O(d^n)$ by Lemmas~\ref{lem:hyperplane-arrangement-bound},~\ref{lem:hyperplane-arrangement}.

	For each maximal cell $P \in \mathcal P_{\max}$, 
    let $I_P := \{ i \in [m] \mid \intr(C_i) \cap P \cap \Z^n \neq \emptyset\}$.
We compute this using \Cref{lem:interior-intersects-lattice} to construct $P'$ such that $\intr(C_i) \cap P \cap \Z^n = \intr(C_i \cap P') \cap \Z^n$ and \Cref{thm:interior-optimization} to detect integer points on the interior of $C_i \cap P'$, and,   moreover we conclude that for each $i \in I_P$, $\intr(C_i) \cap P$ is full-dimensional.
  \medskip

	\noindent \underline{Intermediate step.}
	We claim that 
	\begin{equation}
		\label{eq:polytope-decomposition-integer}
		\Z^n \setminus \bigcup_{i=1}^m \intr(C_i)
		= \bigcup_{P \in \mathcal P_{\max}} \left((P \cap \Z^n) \setminus \bigcup_{i \in I_P} \intr(C_i)\right).
	\end{equation}
	
	To prove this, let \( \x \in \Z^n \). 
	Since \( \{P\}_{P \in \mathcal P_{\max}} \) tiles \( \R^n \), it does so for $\Z^n$ and there exists at least one \( P \in \mathcal P_{\max} \) such that \( x \in P \). 
	
	\textbf{Case 1:} \( \x \notin \bigcup_{i=1}^m \intr(C_i). \)
	
	In this case, the point \( x \) is outside the interiors of all the sets \( C_i \). 
	Therefore, \( x \in \Z^n \setminus (\bigcup_{i=1}^m \intr(C_i))\). 
	For checking the right-hand side (RHS) of~\eqref{eq:polytope-decomposition-integer}, we recall from above that there exists $P \in \mathcal{P}_{\max}$ such that $\x \in P$.
	In this case, $\x \in P \setminus \bigcup_{i=1}^m \intr(C_i) \subseteq P \setminus \bigcup_{i \in I_P} \intr(C_i)$.
	As the RHS is a union of all such terms over $P$, it follows that $\x$ is an element of it.
	
	\textbf{Case 2:} \( \x \in \intr(C_i) \) for $i \in I \subseteq [m]$, $|I| > 0$.
	
	Since $\x \in \bigcup_{i \in I} \intr(C_i) \subseteq \bigcup_{i = 1}^m \intr(C_i)$, $\x$ is not an element of the left-hand side of~\eqref{eq:polytope-decomposition-integer}.
	For testing containment in the RHS, let $P$ be a maximal cell such that $\x \in P$, which exists. 
	Since $\x \in \intr(C_i)$ for $i \in I$, it follows that $P \cap \intr(C_i) \neq \emptyset$ for $i \in I$ and thus $I \subseteq I_P$.
	We conclude that $\x \notin (P \cap \Z^n) \setminus \bigcup_{i \in I_P} \intr(C_i)$.
	Since this holds for all $P$ with $\x \in P$, the point $\x$ is no element of the RHS.
	
	Since both cases hold, we conclude that \eqref{eq:polytope-decomposition-integer} is valid.\\

	\noindent \underline{Restrict to a single cell.}
	We now restrict our attention to some $P\in \mathcal P_{\max}$, and without loss of generality, $I_P = \{1, \dots, m'\}$ for some $m' \leq m$. 
	Since $P$ is defined by a \HC, it holds true that for all $i,\, j \in I_P$, $i \neq j$, we get $(\partial C_i \cap \partial C_j) \cap P \subseteq \partial P$.
	In particular, by noting that $\partial(C_i \cap P)$ is a union of parts of $\partial C_i$ and parts of $\partial P$, we can derive that
	\begin{equation}
		\label{eq:boundary-intersect-on-cell-integer}
		\partial(C_i \cap P) \cap \partial(C_j \cap P) \subseteq \partial P.
	\end{equation}
	
	\textbf{Intersection detection.} 
	Let $C_i^\mathrm{I} := \intr(C_i) \cap P \cap  \Z^n$ for  $i \in [m']$.
	Now, we build a graph $G = ([m'], E)$ with nodes $[m']$ and edges $(i,j) \in E$ whenever $C^\mathrm{I}_i$ and $C^\mathrm{I}_j$ are not weakly separable.
	We detect a weak separation via $\WEAK{\mathcal C}{\Z^n}$ as constructed in~\Cref{lem:separation-integer-hulls} where $\mathcal C = \{P \cap \bigcap_{i \in I'} C_i : P \text{ is a polyhedron and } I' \subseteq I\}$.

	a. If $G$ is connected, then we can compute a connected ordering via a breadth-first-search on $G$ in time $O((m')^2)$; that is, we  compute a permutation of $[m']$ such that $C_i^\mathrm{I} \cap P$ and $C_{i+1}^\mathrm{I} \cap P$ are not weakly separable for $i \in [m' - 1]$.
	Due to the convexity of $C_i$ this also implies that $C_i \cap P$ and $C_{i+1} \cap P$ are not weakly separable and, in other terms, $\bigcup_{j= 1}^i C_j \cap P$ and $C_{i+1} \cap P$ are not weakly separable, $i \in [m' - 1]$.
	Together with~\eqref{eq:boundary-intersect-on-cell-integer}, this allows to apply Corollary~\ref{cor:boundary} for $P$ and its convex subsets $C_i \cap P$, $i \in [m']$.
	This gives that $\bigcup_{i = 1}^{m'} \intr(P \cap C_i) = \intr(P) \cap \bigcup_{i=1}^{m'} \intr(C_i)$ is convex. 
	In particular, $K := P \cap \bigcup_{i = 1}^{m'} \intr(C_i)$ is convex and we can return $(K,\,P,\, [m'])$ as one of our components in the decomposition.

	b. Alternatively, if $G$ is disconnected, then on each connected component, we can make a similar argument that each connected component forms a convex set. 
	Let $K_1, \dots, K_{m''}$ be these convex sets with $2 \leq m'' \leq m' \leq m$.
	Then, for all $v \in [m'']$, there exist  indices $I_v \subseteq [m']$ such that %
	$$
	K_v = P \cap \bigcup_{i \in I_v} \intr(C_i).
	$$
	Without loss of generality, consider $I_1$ and $I_2$.
	Then, for all $i \in I_1$ and $j \in I_2$, the sets $C_i^\mathrm{I} \cap P$ and $C_j^\mathrm{I} \cap P$ are weakly separable, compare above. 
	Equivalently, we can state that $\intr(C_i^\mathrm{I} \cap P) \cap \intr(C_j^\mathrm{I} \cap P) = \emptyset$.
	It follows that 
	$$
	\left( \bigcup_{i \in I_1} \intr(C_i^\mathrm{I} \cap P) \right) \cap \left(\bigcup_{j \in I_2} \intr(C_j^\mathrm{I} \cap P)\right) = \emptyset.
	$$

	Note that these unions contain the integers of the convex sets $K_1$ and $K_2$, respectively.
	Hence, $K_1^\mathrm{I} := K_1 \cap \Z^n$ and $K_2^\mathrm{I} := K_2 \cap \Z^n$ are weakly separable.  We employ $\WEAK{\mathcal C}{\Z^n}$ once more for this.  As $K_i$ is a union of interiors of convex sets, we need an optimization oracle for $K_i$ to employ the lemma.  This is achieved by simply optimizating over each convex set individually using \Cref{thm:interior-optimization} and then comparing the outputs. 
	This approach is valid for any two $K_v$ and $K_w$, $v,\, w \in [m'']$, $v \neq w$.

	\textbf{Separation.} 
	For this, we compute a separating hyperplane $H_{\a^{vw},b_{vw}}$ of $K_v^\mathrm{I}$, $K^\mathrm{I}_w$ for any pair $v,\,w \in [m'']$, $v < w$.
	Note that these are ${ m'' \choose 2}$ many and this is performed via Lemma~\ref{lem:separation-integer-hulls}.
	
	Without loss of generality, we orientate the corresponding halfspaces such that $K_v^\mathrm{I} \subseteq H_{\a^{vw},b_{vw}}^{\leq}$ if $v < w$ and $K_v^\mathrm{I} \subseteq H_{\a^{wv},b_{wv}}^{\geq}$ if $w < v$.
	Now, we set 
	\begin{equation}
		\label{eq:Pv-def-integer}
		P_v := P \cap \bigcap_{v < w} H_{\a^{vw},b_{vw}}^{\leq} \cap \bigcap_{w < v} H_{\a^{wv},b_{wv}}^{\geq},
	\end{equation}
	for all $v \in [m'']$, where $w$ iterates over $[m'']\setminus\{v\}$.
	For the decomposition, we can either restrict each $K_v$ according to the computed hyperplanes above or return $K_v^\mathrm{I}$.
	This gives $(K_v,\, P_v,\, I_v)$ as components in the decomposition. 
	
	As there remain polyhedra in the scheme of the right-side of~\labelcref{eq:Pv-def-integer} when following all the other $v-w$-combinations, we detect all such non-empty ones by means of~\Cref{lem:hyperplane-arrangement}, say $P'_1,\dots,P'_{s'}$.
	We can return them with $K'_1 = \dots = K'_{s'} = \emptyset$ and index sets $I'_1 = \dots = I'_{s'} = \emptyset$.
	
	Since the amount of hyperplanes required for the definition all polyhedra in this paragraph is $m''(m''-1)/2$, as noted above, by~\Cref{lem:hyperplane-arrangement-bound}, we have that $m'' + s' = O(m^{2n})$.
	Since we fixed a cell $P$ and there are at most $O(d^n)$ such, we receive $O((m^2d)^n)$ many pieces. 
	Then we have that 
	\begin{equation}
		\label{eq:P-decomposition-integer}
		(P \cap \Z^n)  \setminus \bigcup_{i \in I_P} \intr(C_i) = \left(\bigcup_{s = 1, \dots ,m''}   (P_s \cap \Z^n) \setminus K_s \right) \cup \left(\bigcup_{s = 1, \dots ,s'}   (P'_s \cap \Z^n) \setminus K'_s\right).
	\end{equation}
	
	Combining \eqref{eq:polytope-decomposition-integer} and \eqref{eq:P-decomposition-integer} completes the result.
\end{proof}

\section{Complexity of Reverse Convex Integer Programming}
\label{sec:complexity-rev-convex}

With the proof ideas from above, we can directly use~\Cref{thm:dd-min} to show that we can solve the feasibility problem over 
$P \setminus \bigcup \intr(C_i)$ in oracle polynomial time. It is natural to subtract open sets in reverse convex programming, 
as this ensures that the resulting feasible region remains closed, preserving the properties needed for efficient computation. 
Furthermore, if $f(x)$ is convex, encoding the reverse convex constraint involves removing the region where $f(x) \geq 0$, 
effectively subtracting the convex set defined by $f(x) < 0$ using strict inequality.

\begin{theorem}[Reverse Convex IP in Fixed-Dimension]
\label{thm:main-IP-complexity}
    Let $P$ be a rational polytope.  Suppose $C_1, \dots, C_m \subseteq \R^n$ are full-dimensional closed convex sets equipped with interior separation oracles $\SEPint{C_i}{\Z^n}$ and have a \HC $\mathcal H$.   
    Then, when the dimension $n$ is fixed, in oracle-polynomial time that is polynomial in $\langle P\rangle, \langle R\rangle, m, \langle \mathcal H\rangle$, we can determine whether the set 
    $$
    (P \setminus (\intr(C_1) \cup \dots \cup \intr(C_m))) \cap \Z^n
    $$ 
    is non-empty.  If so, we can return such a point.
\end{theorem}
\begin{proof}
We first apply~\Cref{thm:main-integer} to obtain a list of polytopes and convex sets $(P'_s, K_s)$ for $s = 1, \dots, t$ with $t \leq O((m^2d)^n)$.  Then, we further intersect each $P'_s$ with $P$, resulting in polyhedra $P_s$.
With choosing $Q_j$ empty in~\Cref{def:division-description}, this gives a \CCPsub.
Since $P$ is a polytope, and hence bounded, the encoding size of the hyperplanes $H_{\a^{vw}, b_{vw}}$ constructed in ~\Cref{thm:main} can be bounded in the encoding size of $P$.  
Hence, the result then follows immediately from~\Cref{thm:dd-min}.
\end{proof}
The proof uses \Cref{thm:main-integer} for the decomposition.  Of course, the practical complexity of this result is dramatically improved if either \Cref{thm:main} or \Cref{thm:ideal-decomposition} can be used.  This depends on what is known about the convex sets.

\section{Characterizations of Intersections of Sets}
\label{sec:characterizations}
In this section, we investigate common types of convex sets and give a sufficient structure for a \HC to exist.

\subsection{Two General Sets by Function Compositions}

\begin{lemma}[General Structure]
\label{lem:general-structure}
Suppose that 
$$
C_0 = \{ \x : f(\x) \leq 0\}, \quad C_1 = \{\x : f_1(\x) \leq 0\},
$$
and 
\begin{equation*}
\label{eq:general_form}
f_1(\x) = g(\x)G(f(\x)) + \prod_{j \in J} h_j(\x),
\end{equation*}
where $f,\, g \colon \R^n \to \R$ are arbitrary functions, $G\colon \R \to \R$, $G(0) = 0$, and  \(h_j\) are affine and non-constant functions for \(j \in J = \{1, \dots, \ell\}\).
Then \(C_0\) and \(C_1\) have a  \HC.
\end{lemma}

\begin{proof}
We show that \(\bd C_0 \cap \bd C_1\) is contained in a collection of hyperplanes.
\begin{align*}
\bd C_0 \cap \bd C_1 &= \{\x : f(\x) = 0,\, f_1(\x) = 0\} \\
&= \{\x : f(\x) = 0,\, g(\x)G(f(\x)) + \prod_{j \in J} h_j(\x) = 0\} \\
&= \{\x : f(\x) = 0,\, \prod_{j \in J} h_j(\x) = 0\} \\
&= \{\x : f(\x) = 0\} \cap \bigcup_{j \in J} \{\x : h_j(\x) = 0\}.
\end{align*}
Since \(h_j\) are affine functions, this last line is indeed a union of hyperplanes.
\end{proof}

This lemma provides one sufficient condition for two sets.  It can be used to ``link" potentially many sets together to create complicated structures to deal with.
We next show how this result appears for quadratic functions.

\subsection{Ellipsoids and Other Quadratic Sets}
Given two ellipsoids, without loss of generality, we can apply an affine transformation so that one ellipsoid is in fact the unit ball, and the other is a general ellipsoid.

In the following, we describe when the intersection of the unit ball with a quadratic set has the property \HC. We then characterize when this quadratic has a positive definite representation, i.e., must be an ellipsoid.

\begin{lemma}
\label{thm:exact-characterization}
    Let $B = \{\x : \|\x\|_2 \leq 1\}$ and $C = \{\x : g(\x) \leq 0\}$ where $g$ is a quadratic polynomial.
Then  $B$ and  $C$ have a \HC if  
\begin{subequations}
\label{eq:g-definitions}
     \begin{equation}
         \label{eq:g-defn}
g(\x) =   \alpha\left(\sum_{i=1}^n x_i^2 - 1\right) + h_1(\x)h_2(\x)
\end{equation}
or
\begin{equation}
\label{eq:g-defn2}
g(\x) =  \alpha\left(\sum_{i=1}^n x_i^2 - 1\right) + h_1(\x),
\end{equation}
\end{subequations}
where  $h_i(\x)$ are  non-trivial affine functions. The hyperplanes are $h_i(\x) = 0$ for either $i=1,2$ in the first case or just for $i=1$ in the second case.
\end{lemma}
\begin{proof}
Suppose one of the formulae for $g$ holds and note that ${\|\x\|_2 = 1 \iff \sum_{i = 1}^n x_i^2 = 1}$. Then
by Lemma~\ref{lem:general-structure}, $B$ and $C$ have \HC with two (possibly equivalent) hyperplanes, $h_1, h_2$.
\end{proof}
We suspect that this theorem is exact, but need some algebraic tools to prove this.
\begin{conjecture}
    \Cref{thm:exact-characterization} holds as an if and only if.
\end{conjecture}

In \cite{BELOTTI20132778}, they study characterizations of quadratic surfaces given the fact that they have the same intersection on two hyperplanes.  
The purpose of this work lies in understanding convexifications of so-called quadrics intersected with linear inequalities.  
Fortunately, they provide somewhat thorough description of the types of surfaces that arise.  They show similar results to \Cref{thm:exact-characterization}.

In order to apply our results, we need the function $g$ to also be convex.  
For this, 
we next do some similar calculations to~\cite{BELOTTI20132778} to discuss when these other functions are convex.

\begin{proposition}[Convex intersections with the ball]
Let $g,h_1, h_2$ be as in  \eqref{eq:g-defn}. In particular,
\begin{align*}
h_1(\x) = \a^\top\x + c_1, \ \ 
h_2(\x) = \b^\top\x + c_2,
\end{align*}
where $\a, \b \not\equiv 0$.
Then $g$ is convex if and only if $ \a^\top\b \geq \|\a\|_2 \|\b\|_2 - 2\alpha$.
\end{proposition}

\begin{proof}
The function \( g(\x) \) is defined as in \eqref{eq:g-defn}.
The Hessian matrix \( \mathbf{H} \) for \( g(\x) \) is given by:
\[
\mathbf{H}_{ij} = 2\alpha \delta_{ij} + a_i b_j + a_j b_i,
\]
where \( \delta_{ij} \) is the Kronecker delta. More compactly, we can express this as:
\[
\mathbf{H} = 2\alpha I + (\a\b^\top + \b\a^\top),
\]
where \( I \) is the identity matrix.

For \( g(\x) \) to be convex, the Hessian \( \mathbf{H} \) must be positive semi-definite (PSD). Since both \( 2\alpha I \) and \( \a\b^\top + \b\a^\top \) are symmetric matrices, \( \mathbf{H} \) is also symmetric. Therefore, \( g(\x) \) is convex if and only if all eigenvalues of \( \mathbf{H} \) are non-negative.

The eigenvalues of \( \mathbf{H} \) are as follows:
\begin{itemize}
\item \( 2\alpha \) with multiplicity \( n - 2 \), corresponding to eigenvectors orthogonal to both \( \a \) and~\( \b \).
\item \( 2\alpha + \a^\top\b + \|\a\|_2 \|\b\|_2 \) with eigenvector proportional to \( \a + \frac{\|\a\|_2}{\|\b\|_2}\b \).
\item \( 2\alpha  
 +\a^\top\b - \|\a\|_2 \|\b\|_2 \) with eigenvector proportional to \( \a - \frac{\|\a\|_2}{\|\b\|_2}\b \).
\end{itemize}

For \( \mathbf{H} \) to be positive semi-definite, the smallest eigenvalue must be non-negative. Therefore, we require:
\[
2\alpha  
 +\a^\top\b - \|\a\|_2 \|\b\|_2\geq 0.
\]
This simplifies to:
\[
\a^\top\b \geq \|\a\|_2 \|\b\|_2 - 2\alpha.
\]

\end{proof}

\begin{corollary}
\label{cor:quadratic-bhc}
Let $B = \{\x : \|\x\|_2 \leq 1\}$ and $C = \{\x : g(\x) \leq 0\}$, where $g$ is a quadratic polynomial of the form as in \eqref{eq:g-definitions}.
If $g$ is of the form~\eqref{eq:g-defn} and  $ \a^\top\b \geq \|\a\|_2 \|\b\|_2 - 2\alpha$ or $g$ is of the form~\eqref{eq:g-defn2}, then $B$ and $C$ have a \HC. 
\end{corollary}
We note that this corollary holds under any invertible affine transformation of the variables as well.
More complicated quadratic intersections also exist.
We defer to~\cite{BELOTTI20132778} for formulas dictating how other quadratic surfaces can intersect hyperplanes with the same intersection.

\section{Conclusions and future work}

In this work, we rely on hyperplane arrangements to create a decomposition.  The approach using hyperplane arrangements imposes an added complexity of $O(d^n)$ where $d$ is the number of hyperplanes in the arrangement.  It is unclear whether this issue can be circumvented to achieve a fixed-parameter tractable complexity result of the form $f(n) \mathrm{poly}(\cdot)$ as can be achieved (typically) from convex integer programming. Here we would want the dependence on the dimension to be multiplicative in the expression and not have exponents related to dimension  on any other parameters.  One approach is to bound the size of the number used in the hyperplane arrangement: this introduces a bound on the number of possible hyperplanes that becomes a function of dimension and this bound.  This could be useful for fixed-parameter tractable type results that address problems with small coefficients.

In future work, we want to address the mixed-integer setting.  As mentioned in the introduction, working with continuous variables over convex sets inherently introduces issues of approximabilitly for even optimizing continuous functions over convex sets.  Thus, additional care needs to be taken in this setting.
Lastly, we conjecture that a version of \Cref{thm:main} holds where we remove closed convex sets instead of interiors, but require that the sets be strictly convex.

\newpage

\printbibliography

\newpage

\appendix

\section{Alternative Proof of Lemma~\ref{lem:boundary}}
\label{sec:alternative-proof}
Here we provide an alternative proof of Lemma~\ref{lem:boundary} that essentially relies on just the Intermediate Value Theorem.  To make this proof precise, however, we also need to work with a specific continuous function called the \emph{Gauge function} (a.k.a., Minkowski Gauge function).

\begin{definition}
\label{defn:gauge}
Let $C \subset \mathbb{R}^n$ be a convex set that is closed, full-dimensional, and contains the origin in its interior. The \textbf{Gauge function} (a.k.a. Minkowski of Minkowski-Gauge function) associated with $C$ is defined as:
\[
\gamma_C(x) = \inf \{\lambda > 0 : x \in \lambda C\}, \quad x \in \mathbb{R}^n.
\]
\end{definition}

We begin with a folklore result about when this function is continuous.  We provide a reference to a recent article with a proof.
\begin{lemma}{\cite[Lemma 3.7]{Cuong2021}}
\label{lem:gauge-continuous}
Let $X$ be a locally convex topological vector space (e.g., $\R^n$).
If $F$ is a convex subset of $X$ with $0 \in \operatorname{int}(F)$, then the Minkowski gauge function is continuous on $X$.
\end{lemma}

\begin{proof}[Alternative proof of Lemma~\ref{lem:boundary} using Intermediate Value Theorem]
Suppose that $C_1$ and $C_2$ are not weakly separable.  Then the intersection $C_1 \cap C_2$ cannot be contained in a hyperplane.  Thus the intersection is full-dimensional and therefore 
there exists $\y \in C_1 \cap C_2$ and $\varepsilon > 0$ such that 
\begin{equation*} 
\label{eq:inner-y}
B(\y, \varepsilon) \subseteq \intr(C_1 \cap C_2),
\end{equation*}
where $B(\y, \varepsilon)$ denotes the open ball around $\y$ with radius $\varepsilon$.

We will prove that $C_1 \cup C_2$ is convex by contradiction.
So, assume $C_1 \cup C_2$ is not convex. 
Then, there exist $\x_1 \in C_1$ and $\x_2 \in C_2$ such that $\z_\lambda = \lambda \x_1 + (1-\lambda)\x_2 \notin C_1 \cup C_2$ for some $\lambda \in (0,1)$.
Without loss of generality, we can assume that $\x_i \in \partial C_i$ for $i = 1, 2$.
Therefore, $\z_\lambda \notin C_1 \cup C_2$ for all $\lambda \in (0, 1)$.

\emph{
\textbf{Claim:} There exists  $\z^* \in \bd C_1 \cap \bd C_2 \cap \relint(T)$ where $T:=\conv(\{\y, \x^1, \x^2\})$.
}\\

Let $\z_\lambda(t) = t \z_\lambda + (1-t) \y$ for $t \in [0,1]$ and $C_i - \y := \{ \x - \y \mid \x \in C_i\}$, for $i = 1,\, 2$.
As $\y \in \intr(C_i)$, the origin is contained in $\intr(C_i - \y)$, and we can define the respective gauge function $\gamma_{C_i - \y}(\x)$, for $i = 1,\,2$.

Now, let 
$$
f_{i}(\lambda) = \frac{1}{\gamma_{C_i - \y}(\z_{\lambda} - \y)},
$$
for $\lambda \in [0, 1]$.
Note that $\gamma_{C_i - \y}$ continuous due to \Cref{lem:gauge-continuous} and $\gamma_{C_i - \y}(\z_\lambda - \y) > 0$, as $\z_\lambda - \y \neq \textbf{0}$ for such $\lambda$.
Therefore, function $f_i(\lambda)$ is well-defined, continuous and strictly positive.

For now, consider an arbitrary but fixed $\lambda \in [0, 1]$.
Say, 
$$
\lambda_i := 1/f_i(\lambda) =  \gamma_{C_i - \y}(\z_\lambda -\y) = \inf\{\mu > 0 : \z_\lambda - \y \in \mu (C_i - \y)\}.
$$
Then, 
$$
\z_\lambda(f_i(\lambda)) = \z_\lambda(1/\lambda_i) = (1/\lambda_i)\z_\lambda + (1 - (1/\lambda_i)) \y = (1/\lambda_i) (\z_\lambda + (\lambda_i - 1) \y).
$$
Together, we derive by the positive homogeneity of $\gamma_{C_i - \y}$ that
\begin{align*}
    & \z_\lambda - \y \in \bd( \lambda_i(C_i - \y))\\
    \iff & 1/\lambda_i(\z_\lambda - \y) \in \bd(C_i - \y) \\
    \iff & 1/\lambda_i(\z_\lambda + (\lambda_i - 1)\y) \in \bd C_i \\
    \iff & \z_\lambda(f_i(\lambda)) \in \bd C_i.
\end{align*}

Observe that $f_1(1) = 1$ and $f_2(0) = 1$ since $\x_1 \in \bd C_1$ and $\x_2 \in \bd C_2$, respectively.
In addition, since $\x_2 \notin C_1$ and $\x_1 \notin C_2$, it is $\gamma_{C_1 - \y}(\x_2 - \y) > 1$ and $\gamma_{C_2 - \y}(\x_1 - \y) > 1$, thus $f_1(0) < 1$, $f_2(1) < 1$, respectively.

Now, define $g(\lambda) = f_1(\lambda) - f_2(\lambda)$.
From above we know that $g$ is continuous and $g(0) < 0$, as well as $g(1) > 0$.  
Thus by the intermediate value theorem, there exists a $\bar \lambda \in (0,1)$ such that $g(\bar \lambda) = 0$.    In particular, defining $\z^* := \z_{\bar \lambda}(f_i(\bar \lambda))$ we have $\z^* \in \partial C_1 \cap \partial C_2$. 

As $0 < \bar \lambda < 1$, one can verify that $\z^* \in \relint(T)$.
Indeed, $\z_{\bar\lambda}$ is a true convex combination of the vertices $\x_1$, $\x_2$, $\y$ of $T$.
This completes the claim.

Finally, with $\z^* \in \relint(T)$ and  $\relint(T) \subseteq \relint(K)$ we derive that $\z^* \notin \bd K$.
This is a contradiction, since $\z^* \in \bd C_1 \cap \bd C_2 \subseteq \bd K$.

By assumption, we conclude that $C_1 \cup C_2$ is convex if $C_1$ and $C_2$ are not weakly separable.
\end{proof}

\section{Integer Programming with Queries only on Integers}
\label{sec:complexity-integer-queries}
We follow Basu~\cite{Basu2022} to present a self-contained version of the proof for \Cref{thm:convex-integer}.  

Let \( \|\x\|_2 = \sqrt{\x^\top \x} \) denote the standard Euclidean norm. For a symmetric positive definite matrix \( A \in \mathbb{R}^{n \times n} \), the \emph{ellipsoidal norm} induced by \( A \) is defined as
\[
\|\x\|_A := \sqrt{\x^\top A \x}.
\]

We define an ellipsoid centered at \( c \in \mathbb{R}^n \) with shape matrix \( A \succ \ve 0 \) as
\[
E := \c + \mathcal{E}(A) = \{ \x \in \mathbb{R}^n : \|\x - \c\|_A \leq 1 \}.
\]
Its volume is given by
\[
\vol(E) = \frac{\pi^{n/2}}{\Gamma(n/2 + 1)} \cdot \frac{1}{\sqrt{\det(A)}}.
\]

We assume that the convex set \( C \subseteq \mathbb{R}^n \) satisfies \( C \subseteq E \subseteq B_2(\ve 0, R) \), i.e., it lies within an ellipsoid contained in a Euclidean ball of radius \( R \).

We have access to a \emph{separation oracle} for \( C \)
and emphasize that the algorithm only queries the oracle at integer points.

The \emph{Closest Vector Problem (CVP)} in ellipsoidal norm \( \| \cdot \|_A \) is defined as: given \( c \in \mathbb{R}^n \), find
\[
\x^* = \arg\min_{\x \in \mathbb{Z}^n} \|\x - \c\|_A.
\]

The \emph{Shortest Vector Problem (SVP)} in the dual norm \( \| \cdot \|_{A^{-1}} \) is to find a nonzero integer vector \( \w \in \mathbb{Z}^n \setminus \{\ve 0\} \) minimizing \( \|\w\|_{A^{-1}} \). This vector is used to bound the \emph{lattice width} of the ellipsoid.

For a direction \( w \in \mathbb{Z}^n \setminus \{0\} \), the lattice width of a convex set \( K \subseteq \mathbb{R}^n \) in direction \( \w \) is defined as
\[
\mathrm{width}_\w(K) := \max_{\x \in K} \langle \w, \x \rangle - \min_{\x \in K} \langle \w, \x \rangle.
\]
If \( \mathrm{width}_w(K) < 1 \), then \( K \cap \mathbb{Z}^n \) lies in the hyperplane \( \{\x \in \mathbb{R}^n : \langle \w, \x \rangle = k\} \) for some integer \( k \), by Basu~\cite[Lemma 5.5]{Basu2022}.

Throughout, we assume that CVP and SVP are solved using the deterministic \( 2^{O(n)} \)-time algorithm of Micciancio and Voulgaris~\cite{MicciancioVoulgaris2010}, and that ellipsoids and separating hyperplanes are rounded appropriately using the techniques of Frank and Tardos~\cite{FrankTardos1987} to maintain polynomially bounded bit-size.

\begin{proof}[Proof of \Cref{thm:convex-integer} following \cite{Basu2022}]
Let \(C \subseteq B_\infty(\ve 0,R)\) be a convex set given by the integer separation oracle
\(\SEP{C}{\Z^n}\). We design a recursive ellipsoid-based algorithm that decides whether 
\(C \cap \Z^n = \emptyset\) using only integer separation queries. 
The algorithm is built as a recursive function $\IE_n(C,E)\), where 
\(E = \mathbf{c} + \mathcal{E}(A) \subseteq \R^n\) is an ellipsoid containing \(C\). 
The outcome of \(\IE_n(C,E)\) is defined recursively as one of four possibilities:
\[
\begin{cases}
\texttt{(Found)}\ \x & \text{for some } \x \in C \cap \Z^n, \\[6pt]
\texttt{(Rounding)}\ \IE_n(C,E') & \text{where } C \subseteq E' \text{ and } 
   \vol(E') \leq \rho_n \cdot \vol(E), \\[10pt]
\texttt{(Subcases)}\ \{\IE_{n-1}(C_k,E_k)\}_{k=\ell}^u    &\text{where } 
   \w \in \Z^n \setminus \{\ve 0\} \text{ is a flatness direction and} \\[6pt]
&   C_k = \proj_{n-1}(C \cap \{\x : \langle \w,\x \rangle = k\}), \\[6pt]
   &   E_k = \proj_{n-1}(E \cap \{\x : \langle \w,\x \rangle = k\}), \\[6pt]
\texttt{(Empty)} & \text{assertion that } C \cap \Z^n = \emptyset, \\
   & \text{triggered when all recursive branches are exhausted.}
\end{cases}
\]

where 
\(
\rho_n = \exp\!\Big(-\tfrac{1}{5n(n+1)^2}\Big)
\)
and where $\proj_{n-1}$ projects into a new space with $n-1$ many variables via a change of coordinates, typically computed using a Hermite normal form computation.  

We now explain the steps and analysis of the procedure.

\medskip

\textbf{Step 1: Volume Check.}  
If \(\vol(E) < 1/n!\), then by \cite[Lemma 5.5]{Basu2022}, 
all integer points in \(C\) lie on a common affine hyperplane. 
In this case we output \texttt{(Subcases)} by returning a flatness direction \(\w\) 
and a single slice \([\ell,u]=[\delta,\delta]\) where this output is computed using \textsc{FlatnessRecursion}\((C,E)\), as defined below.

\textbf{Step 2: Closest Integer Point.}  
Otherwise (\(\vol(E) \ge 1/n!\)), compute the closest integer point to the ellipsoid center \(\mathbf{c}\) in the ellipsoidal norm:
\[
\x^\star := \arg\min_{\x \in \Z^n} \|\x - \mathbf{c}\|_A.
\]

\textbf{Case A (Close Integer Point).}  
If \(\|\x^\star - \mathbf{c}\|_A \le 1/(n+1)\), then \(\x^\star \in E\).  
\begin{itemize}
    \item If \(\SEP{C}{\Z^n}(\x^\star) \) conforms $\x^* \in C$, return outcome \texttt{(Found)} with \(\x^\star \in C \cap \Z^n\).
    \item Otherwise the oracle provides a separating hyperplane \(\a^\top \x \le \beta\) valid for \(C\) but violated by \(\x^\star\). 
    By \cite[Lemma 5.6]{Basu2022}, construct an ellipsoid \(E' \supseteq C\) as the minimum volume ellipsoid enclosing $E \cap \{\x : \a^\top \x \le \beta\}$ which has 
    \(\vol(E') \leq \rho_n \cdot \vol(E)\), and return outcome \texttt{(Rounding)} with \(E'\). 
\end{itemize}

\textbf{Case B (Far Integer Point).}  
If \(\|\x^\star - \mathbf{c}\|_A > 1/(n+1)\), invoke \textsc{FlatnessRecursion}\((C,E)\).

\medskip
\textbf{Subroutine: \textsc{FlatnessRecursion}\((C,E)\).}
\begin{enumerate}
    \item Solve the shortest vector problem (SVP) in the dual ellipsoidal norm 
    \(\|\z\|_{A^{-1}} = \sqrt{\z^\top A \z}\) to obtain a flatness direction
    \(\w \in \Z^n \setminus \{\ve 0\}\) minimizing
    \[
    \text{width}_\w(E) = \max_{\x \in E} \langle \w,\x \rangle - \min_{\x \in E} \langle \w,\x \rangle.
    \]
    \item Let \([\ell,u] \subseteq \Z\) be the interval of integer slices intersecting \(E\):
    \[
    C \cap \Z^n \;\subseteq\; \bigcup_{k=\ell}^u \{\x \in \Z^n : \langle \w,\x \rangle = k\}.
    \]
    \item \textbf{Subcase B1 (Very Flat).} If \(\text{width}_\w(E) < 1\), then by \cite[Lemma 5.5]{Basu2022} all integer points lie in a single hyperplane 
    \(\langle \w,\x\rangle = \delta\). We return outcome \texttt{(Subcases)} with \((\w,[\delta,\delta])\).
    \item \textbf{Subcase B2 (Flat but not Collapsed).} Otherwise, \cite[Lemma 5.8]{Basu2022}\ ensures
    \[
    u - \ell + 1 \;\leq\; n(n+1).
    \]
    We return outcome \texttt{(Subcases)} with \((\w,[\ell,u])\).
\end{enumerate}

\medskip
\textbf{Correctness and Termination.}  
Each outcome guarantees progress:
\begin{itemize}
    \item \texttt{(Found)} halts with a feasible integer point.
    \item \texttt{(Rounding)} strictly reduces volume by a factor \(\rho_n < 1\).
    \item \texttt{(Subcases)} reduces the dimension by recursion on at most \(n(n+1)\) slices.
    \item \texttt{(Empty)} is triggered only when all recursive branches in dimension \(<n\) 
    are exhausted, certifying \(C \cap \Z^n = \emptyset\).
\end{itemize}

\textbf{Runtime Analysis.}  
There are at most \(O(n^4 \log R)\) consecutive \texttt{(Rounding)} steps before flatness is forced with $\vol(E) < \tfrac{1}{n!}$.  
Each \texttt{(Subcases)} call branches into at most \(n(n+1)\) instances in dimension \(n-1\).  
Thus, if \(T(n)\) denotes the number of recursive calls in dimension \(n\), then
\[
T(n) \;\leq\; O(n^4 \log R) \;+\; n(n+1)\,T(n-1), 
\qquad T(1) = O(\log R).
\]
By induction,
\[
T(n) \;\leq\; 2^{O(n\log n)} \cdot \log R.
\]
Each call to \textsc{IE-Step} may solve either a CVP or SVP, both computable in time \(2^{O(n)}\) \cite{MicciancioVoulgaris2010}, and may invoke the separation oracle $\SEP{C}{\Z^n}$.  
Hence the total runtime is
\[
2^{O(n\log n)} \cdot \operatorname{poly}\big(\langle R \rangle, \Phi_{\SEP{C}{\Z^n}}(n\langle R \rangle)\big).
\]

\textbf{Optimization.}  
Optimization \(\max/\min \{\c^\top \x : \x \in C \cap \Z^n\}\) follows by binary search on the objective, adding the constraint \(\c^\top \x \ge \gamma\) (or \(\le \gamma\)) at each step.  
This preserves the oracle model and does not change the order of complexity.
\end{proof}

\textbf{Bit Complexity Remark.} At each step of the algorithm, we must ensure that the ellipsoid parameters (center \( \c \) and shape matrix \( A \)) and the separating hyperplanes remain polynomially bounded in size. This can be guaranteed by standard rounding techniques and the framework developed by Frank and Tardos~\cite{FrankTardos1987}, which ensure that all quantities can be encoded with bit-length polynomial in the input size (i.e., in \( \log R \) and the encoding length of \( \SEP{C}{D} \)). Therefore, the storage complexity of the algorithm is also polynomially bounded throughout the recursion.

\begin{remark}[Comparison to Ellipsoid Rounding in Lenstra-Type Algorithms]
Traditional Lenstra-type algorithms for integer programming over convex sets, such as those in Gr{\"o}tschel, Lov{\'a}sz, and Schrijver~\cite{Grtschel1993}, begin with an \emph{ellipsoid rounding} step. This step constructs a pair of concentric ellipsoids \( E^- \subseteq C \subseteq E^+ \) such that the ratio of their volumes is bounded, and this is achieved by querying the separation oracle \( \SEP{C}{D} \) at rational (non-integer) points. This containment allows the algorithm to obtain a bounded-volume representation of \( C \) before performing any operations over the integer lattice.

However, in our setting, we assume that the separation oracle \( \SEP{C}{D} \) can only be queried at \emph{integer} points. In this regime, classic ellipsoid rounding is not directly applicable. To manage this, alternative techniques such as \emph{centerpoints} from Oertel~\cite{OERTEL2014424} and Oertel and Basu~\cite{oertel} can be used to simulate the effect of a volumetric center based on integer oracle access.

In this work, we follow the approach of Basu et al.~\cite{Basu2022}, which bypasses explicit centerpoint computation by instead working with a containment ellipsoid \( E \supseteq C \). Rather than asking for a shallow cut first, the algorithm begins by solving a full \emph{Closest Vector Problem (CVP)} in the ellipsoidal norm to find the nearest integer point to the center of \( E \), and only then queries the oracle.

While this appears more computationally demanding—since we must solve CVP at every step before getting a cut—it avoids the need for full ellipsoid rounding. Moreover, the overhead of additional CVP calls (each solvable in \( 2^{O(n)} \) time via Micciancio and Voulgaris~\cite{MicciancioVoulgaris2010}) is ultimately dominated by the number of recursive branches in the full algorithm (e.g., \( n(n+1) \) branches in flatness recursion). Thus, the overall asymptotic complexity remains \( 2^{O(n \log n)} \cdot \operatorname{poly}(\langle R\rangle) \), and this structure supports our restriction to integer-only oracle access.
\end{remark}

\end{document}

%% file: AN1.tex
 \begin{tikzpicture}[scale = 0.7]
    \draw[very thin, gray!30, step=1 cm](-0.3,-0.3) grid (6.9,3.9);

    \fill [blue!25, opacity = 0.7, domain=sqrt(2):5, samples = 100,variable=\x]
      (1.414, 0)
      -- plot ({\x}, {sqrt(0.2)*sqrt(\x*\x-2)+sqrt(0.2)})
      -- (5, 0)
      -- cycle;

     \fill [green!25, opacity = 0.7, domain=sqrt(2):5, samples = 100,variable=\x]
      (2, 0)
      -- plot ({\x}, {sqrt(0.2)*sqrt(\x*\x-2)})
      -- (5, 0)
      -- cycle;

    \draw [thick] [->] (0,0)--(7,0) node[right, below] {$x$};
     \foreach \x in {0,...,6}
       \draw[xshift=\x cm, thick] (0pt,-1pt)--(0pt,1pt) node[below] {$\x$};

    \draw [thick] [->] (0,0)--(0,4) node[above, left] {$y$};
     \foreach \y in {0,...,3}
       \draw[yshift=\y cm, thick] (-1pt,0pt)--(1pt,0pt) node[left] {$\y$};

    \draw [domain=sqrt(2):5, variable=\x, samples = 200]
      plot ({\x}, {sqrt(0.2)*sqrt(\x*\x-2)}) node[right] at (3.2,0.5) {$C$};
          \draw [domain=sqrt(2):5, variable=\x, samples = 100]
      plot ({\x}, {sqrt(0.2)*sqrt(\x*\x-2) + sqrt(0.2)}) node[right] at (3.2,2.5) {$K$};
\draw (5,2.591)--(5,0);
\draw (1.414,0.4472)--(1.414,0);
  \end{tikzpicture}

%% file: pell-pic.tex
 \begin{tikzpicture}[scale = 0.7]
    \draw[very thin, gray!30, step=1 cm](-0.3,-0.3) grid (6.9,3.9);
    
    \fill [blue!25, opacity = 0.7, domain=0:7, samples = 100,variable=\x]
      (0, 4)
      -- plot ({\x}, {sqrt(0.2)*sqrt(\x*\x+2)})
      -- (7, 4)
      -- cycle;

    \fill [green!25, opacity = 0.7, domain=sqrt(2):7, samples = 100,variable=\x]
      (2, 0)
      -- plot ({\x}, {sqrt(0.2)*sqrt(\x*\x-2)})
      -- (7, 0)
      -- cycle;

    \draw [thick] [->] (0,0)--(7,0) node[right, below] {$x$};
     \foreach \x in {0,...,6}
       \draw[xshift=\x cm, thick] (0pt,-1pt)--(0pt,1pt) node[below] {$\x$};

    \draw [thick] [->] (0,0)--(0,4) node[above, left] {$y$};
     \foreach \y in {0,...,3}
       \draw[yshift=\y cm, thick] (-1pt,0pt)--(1pt,0pt) node[left] {$\y$};

    \draw [domain=sqrt(2):7, variable=\x, samples = 200]
      plot ({\x}, {sqrt(0.2)*sqrt(\x*\x-2)}) node[right] at (3.2,0.5) {$C_2$};
          \draw [domain=0:7, variable=\x, samples = 100]
      plot ({\x}, {sqrt(0.2)*sqrt(\x*\x+2)}) node[right] at (3.2,2.5) {$C_1$};

  \end{tikzpicture}

%% file: proof-picture.tex
\begin{figure}
\begin{center}
\begin{tikzpicture}[scale = 0.5]

\definecolor{greenish}{RGB}{144,238,144}
\definecolor{blueish}{RGB}{173,216,230}

\fill[greenish, opacity=0.5] (0,2) ellipse (4cm and 1.5cm);

\fill[blueish, opacity=0.5] (0,0) ellipse (4cm and 1.5cm);

\draw[fill=black] (4,0) circle (2pt) node[anchor= west] {$x_2$};
\draw[fill=black] (4,2) circle (2pt) node[anchor= west] {$x_1$};

\draw[fill=black] (1,1) circle (2pt) node[anchor= east] {$y$};

\draw (1,1) -- (4,0) -- (4,2) -- (1,1);

\draw[fill=purple] (3,1) circle (2pt) node[anchor= east] {$z^*$};
\end{tikzpicture}
\end{center}

\caption{Key ideas in the proof of Lemma~\ref{lem:boundary}.}
\label{fig:proof-lem-boundary}
\end{figure}

%% file: heart.tex
\begin{tikzpicture}[scale = 1.2]
    \fill[blue, opacity = 0.2] (0,0) -- (1,0) -- (1,1.5) arc[start angle=0, end angle=180, radius=0.5] -- (0,0);
    \fill[green, opacity=0.2] (0,0)--(0,1) -- (1.5,1) arc[start angle=90, end angle=-90, radius=0.5] -- (0,0);

    \draw[dotted] (-0.5,0) -- (2.5,0); %
    \draw[dashed](0,-0.5) -- (0,2.5); %
    \draw(2.5,-0.5) -- (-0.5,2.5); %
\end{tikzpicture}

%% file: heart1.tex
\begin{tikzpicture}[scale = 1.2]
    \fill[gray, opacity = 1] (0,0) -- (1,0) -- (1,1.5) arc[start angle=0, end angle=180, radius=0.5] -- (0,0);
    \fill[gray, opacity=1] (0,0)--(0,1) -- (1.5,1) arc[start angle=90, end angle=-90, radius=0.5]  -- (0,0);

        \fill[white, opacity=1] (0,2)--(0,3) -- (3,3) -- (3,0) -- (2,0)--cycle;
    
    \draw[dotted] (0,0) -- (2,0); %
    \draw[dashed](0,0) -- (0,2); %
    \draw(2,0) -- (0,2); %
\end{tikzpicture}

%% file: heart2.tex
\begin{tikzpicture}[scale = 1.2]
    \fill[blue, opacity = 0.2] (0,0) -- (1,0) -- (1,1.5) arc[start angle=0, end angle=180, radius=0.5] -- (0,0);
    \fill[green, opacity=0.2] (0,0)--(0,1) -- (1.5,1) arc[start angle=90, end angle=-90, radius=0.5] -- (0,0);

        \fill[white, opacity=1] (0,2)--(0,0) -- (2,0)--cycle;
    
\draw[dotted] (2,0) -- (2.5,0); %
    \draw[dashed](0,2) -- (0,2.5); %
    \draw(2,0) -- (0,2); %
        \draw[orange](1,1) -- (2,2); %
\end{tikzpicture}

%% file: heart2a.tex
\begin{tikzpicture}[scale = 1.2]
    \fill[blue, opacity = 0.2] (0,0) -- (1,0) -- (1,1.5) arc[start angle=0, end angle=180, radius=0.5] -- (0,0);

        \fill[white, opacity=1] (0,2)--(0,0) -- (2,0)--cycle;
    
    \draw[dashed](0,2) -- (0,2.5); %
    \draw(1,1) -- (0,2); %
        \draw[orange](1,1) -- (2,2); %
\end{tikzpicture}

%% file: heart2b.tex
\begin{tikzpicture}[scale = 1.2]
    \fill[green, opacity=0.2] (1,1) -- (1.5,1) arc[start angle=90, end angle=-90, radius=0.5]  -- (1,1);

        \fill[white, opacity=1] (1,1)--(1,0) -- (2,0)--cycle;
    
    \draw[dotted] (2,0) -- (2.5,0); %
    \draw(1,1) -- (2,0); %
        \draw[orange](1,1) -- (2,2); %
\end{tikzpicture}